\documentclass[11pt, reqno]{amsart}
\usepackage[dvipsnames]{xcolor}
\usepackage{subfiles}
\usepackage[utf8]{inputenc}
\usepackage{amscd, amsmath, amssymb, amsthm, yfonts, mathrsfs, hhline, tikz, dsfont, tikz-cd, mathabx, esvect, accents, dsfont}
\usepackage{wasysym}
\usetikzlibrary{arrows, decorations.markings, decorations.pathreplacing, patterns, calc, intersections, through}
\tikzset{>=latex}
\tikzcdset{arrow style=math font}
\usetikzlibrary{decorations.pathmorphing, intersections, shapes}

\usepackage[bookmarks=true, bookmarksopen=true,%
bookmarksdepth=3,bookmarksopenlevel=2,%
colorlinks=true,%
linkcolor=blue,%
citecolor=blue,%
filecolor=blue,%
menucolor=blue,%
urlcolor=blue]{hyperref}

\usepackage[centering,margin=1.2in]{geometry}
\usepackage[matrix,arrow,curve]{xy}
\global\binoppenalty=10000
\usepackage{todonotes,setspace} 
\usepackage{verbatim}
\usepackage{graphicx}
\usepackage[capitalise]{cleveref}
\DeclareMathAlphabet{\mathpzc}{OT1}{pzc}{m}{it}

\newtheorem{theorem}{Theorem}[section]

\newtheorem{lemma}[theorem]{Lemma}

\newtheorem{prop}[theorem]{Proposition}
\newtheorem{cor}[theorem]{Corollary}
\newtheorem{conjecture}[theorem]{Conjecture}

\theoremstyle{definition}
\newtheorem{defn}[theorem]{Definition}
\newtheorem{remark}[theorem]{Remark}
\newtheorem{example}[theorem]{Example}

\numberwithin{equation}{section}

\def\beq{\begin{equation}}
\def\eeq{\end{equation}}

\newcommand{\longra}{\longrightarrow}

\newcommand{\hr}[1]{\left(#1\right)} 
\newcommand{\hm}[1]{\left|#1\right|} 
\newcommand{\ha}[1]{\left\langle#1\right\rangle} 
\newcommand{\hs}[1]{\left[#1\right]} 
\newcommand{\hc}[1]{\left\{#1\right\}} 

\newcommand{\brr}[1]{(\hspace{-1.5pt}(#1)\hspace{-1.5pt})}
\newcommand{\bss}[1]{[\hspace{-1pt}[#1]\hspace{-1pt}]}

\def\le{\leqslant}
\def\ge{\geqslant}

\def\Ac{\mathcal A}
\def\Ad{\operatorname{Ad}}

\def\Acr{\mathscr A}
\def\Abb{\mathbb{A}}

\def\b{\mathfrak b}

\def\bs{\boldsymbol}

\def\C{\mathbb C}

\def\dim{\operatorname{dim}}

\def\Dc{\mathcal D}
\def\Dcr{\mathscr D}

\def\dimv{d}
\def\Dres{\Dc_{\mathrm{res}}}

\def\eff{\xi}
\def\eps{\varepsilon}

\def\End{\operatorname{End}}

\def\Fcr{\mathscr F}
\def\Frac{\operatorname{Frac}}

\def\gl{\mathfrak{gl}}
\def\rg{\mathrm{gauge}}

\def\Gr{\mathrm{Gr}}

\def\Hom{\operatorname{Hom}}

\def\i{\mathbf i}

\def\id{\mathrm{id}}

\def\Kc{\mathcal K}

\def\la{\lambda}

\newcommand{\lra}{\longrightarrow}

\def\loc{\mathrm{loc}}
\def\longra{\longrightarrow}
\def\La{\Lambda}
\def\Lbb{\mathbb L}

\def\mut{\mathrm{mut}}

\def\Nc{\mathcal N}

\def\Oc{\mathcal O}

\def\pt{\mathrm{pt}}

\def\Prm{\mathrm P}

\def\Q{\mathbb Q}
\def\Qc{\mathcal Q}

\def\Qop{\mathbf{Q}}

\def\res{\operatorname{res}}

\def\rg{\mathrm{g}}
\def\rf{\mathrm{f}}

\def\Rc{\mathcal R}
\def\Rcr{\mathscr R}

\def\sgn{\operatorname{sgn}}

\def\Sc{\mathcal S}
\def\St{\mathrm{St}}

\def\supp{\mathrm{supp}}

\def\tgt{\mathfrak{t}}

\def\Tbb{\mathbb T}
\def\Tc{\mathcal T}

\def\Vcr{\mathscr V}

\def\Wbb{\mathbb{W}}
\def\Wc{\mathcal W}

\def\Xc{\mathcal X}

\def\Z{\mathbb Z}

\def\quiver{\Pi}

\title[$K$-theoretic Coulomb branches and cluster varieties]{$K$-theoretic Coulomb branches of quiver gauge theories and cluster varieties}

\author[Gus Schrader]{Gus Schrader}
\author[Alexander Shapiro]{Alexander Shapiro}

\begin{document}

\maketitle

\begin{abstract}
For $\Gamma$ a quiver without 1-cycles, we show that the Braverman--Finkelberg--Najakima quantized $K$-theoretic Coulomb branch algebra $\Acr_\Gamma$ of the corresponding quiver gauge theory is isomorphic to the quantized universally Laurent algebra (upper cluster $\Xc$-algebra) associated to an explicit initial seed $\quiver_\Gamma$. We also show that $\quiver_\Gamma$ admits a cluster Donaldson--Thomas transformation as defined by Keller.
\end{abstract}

\setcounter{tocdepth}{1}

\section{Introduction}
In \cite{Nak16, BFN18}, Braverman, Finkelberg, and Nakajima propose a mathematical definition of the \emph{Coulomb branch} of a $3d$ $\Nc=4$ supersymmetric gauge theory of cotangent type. Given a complex reductive group $G$ and its complex representation $N$, they introduce a moduli space $\Rcr = \Rcr_{G,N}$, which coincides with the affine Grassmannian $\Gr_G=G\brr{z}/G\bss{z}$ in the case $N=0$. It is shown in~\cite{BFN18} that there is a well-defined $G\bss{z}$-equivariant Borel--Moore homology group $H_*^{G\bss{z}}(\Rcr)$ of the moduli space $\Rcr$, and that $H_*^{G\bss{z}}(\Rcr)$ carries a convolution product which equips it with a commutative ring structure. The Coulomb branch $\mathscr{M}_C = \mathscr{M}_C(G,N)$ of the theory is defined to be the spectrum of this ring:
$$
\mathscr{M}_C = \mathrm{Spec}\hr{H_*^{G\bss{z}}(\Rcr)}.
$$
As a byproduct of the above definition, one obtains a \emph{quantized Coulomb branch}.  The latter is a non-commutative convolution algebra
$$
\mathscr{A}^\hbar_{G,N} = H_*^{G\bss{z} \rtimes \C^\times}(\Rcr),
$$
where the torus $\C^\times$ acts by the \emph{loop rotation} scaling the variable $z$ in $\C\bss{z}$, so that its equivariant parameter $\hbar$ becomes the quantization parameter. There is a variant of the construction above obtained using equivariant algebraic $K$-theory in place of equivariant homology, producing a noncommutative algebra 
$$
\Acr^q_{G,N} =  K_{G\bss{z} \rtimes \C^\times}(\Rcr),
$$
where now we write $q$ for the equivariant parameter for the loop rotation torus.
Physically, the $K$-theoretic Coulomb branches of $3d$ $\Nc=4$ gauge theories can be understood as the Coulomb branches of the $4d$ $\Nc=2$ theories compactified on a circle. It is the $K$-theoretic Coulomb branches and their quantization that we study in this article.

Again from the physical point of view, the work of Gaiotto, Moore, and Neitzke (see \cite{GMN10,GMN13}, as well as~\cite{Nei14} for an overview) suggests that the $S^1$-compactified $4d$ $\Nc=2$ Coulomb branches should be equipped with a collection of cluster-type coordinate systems. As explained in~\cite{Nei14}, this expectation has been made precise in the case of the theories of `class $S$', where the cluster coordinates turn out to coincide with those introduced by Fock and Goncharov in~\cite{FG06} to parametrize moduli spaces of framed $PGL_n$ local systems on the marked surface $S$. It has been conjectured by Gaiotto that a corresponding cluster structure also exists in the case of the supersymmetric gauge theories appearing in the work of Braverman, Finkelberg, and Nakajima.
 
The main result of this article is such a description in the context of an extensively studied class of the $4d$ $\mathcal{N}=2$ theories fitting into the paradigm of Braverman, Finkelberg, and Nakajima --- namely, the quiver gauge theories (for quivers without 1-cycles). In this case, both the group $G = G_\Gamma$ and the representation $N = N_\Gamma$ are encoded by the data of a representation of a (framed) quiver~$\Gamma$, and we write $\Acr_\Gamma^q = \Acr_{G,N}^q$. On the mathematical side, many interesting non-commutative algebras including the quantum group $U_q(\gl_n)$ and quantized Slodowy slices in Dynkin type $A$, are known to arise as quantized Coulomb branch rings $\Acr^q_{\Gamma}$ for appropriate quiver theories. In Theorem~\ref{thm:main-iso}, we prove that $\Acr^q_{\Gamma}$ is isomorphic to the quantized universally Laurent algebra $\Lbb_{\quiver_\Gamma}^q$ (also known as the upper cluster $\Xc$-algebra) associated to an explicit initial cluster seed $\quiver_\Gamma$.\footnote{In the initial version of this paper, we identified the seed $\quiver_\Gamma$ and showed that certain subalgebra of the quantized Coulomb branch $\Acr^q_{\Gamma}$ embeds into the universally Laurent algebra $\Lbb^q_{\quiver_\Gamma}$.} In Section~\ref{sec:DT} we also prove that $\quiver_\Gamma$ admits a \emph{cluster Donaldson--Thomas transformation} in the sense of~\cite{Kel17} and present a formula for the corresponding element of the cluster modular group. Finally, when the quiver $\Gamma$ is unframed and has trivial first Betti number, the cluster $\mathcal{X}$-seed $\quiver_\Gamma$ admits an extension to a compatible pair so that $\Acr^q_{\Gamma}$ contains well-defined cluster $\Ac$-variables. In Section~\ref{sec:avars}, we compute the $\Ac$-variables for our initial seed $\quiver_\Gamma$ and show that they consist of \emph{dressed minuscule monopole operators}, i.e.\ classes of sheaves supported on the preimages of closed orbits in the affine Grassmannian.

The proof of Theorem~\ref{thm:main-iso} exploits a description (see Proposition~\ref{prop:conj-crit}) of the quantized Coulomb branch $\Acr_\Gamma^q$ as a subring in the Coulomb branch algebra for the \emph{pure} gauge theory with the same gauge group $G_\Gamma$ and representation $N=0$. The latter Coulomb branch has been shown to admit a cluster structure by Cautis and Williams in~\cite{CW18}. The cluster structure for $\Acr_\Gamma^q$ is obtained by adding extra mutable directions to those from the initial seed for $\Acr_{G_\Gamma,0}$, while removing some of the frozen ones. At the heart of our argument are the relations intertwining the action of certain sequences of mutations involving these new matter directions with the operators appearing in our description of  $\Acr_\Gamma$ from Proposition~\ref{prop:conj-crit}. The mutation sequences in question are the \emph{bi-fundamental Baxter operators}, which also appear in~\cite{SS25} in the context of cutting and gluing in quantum higher Teichm\"uller theory. It is the same Baxter operators which serve as the building block in our description of the cluster Donaldson--Thomas transformation for $\Acr_\Gamma^q$.

We conclude the introduction by mentioning the work~\cite{CW23} of Cautis and Williams, in which they construct a finite-length $t$-structure on the derived category of $G\bss{z}$-equivariant coherent sheaves on the moduli space $\Rcr$, leading to to a canonical basis in the $K$-theoretic Coulomb branch algebra. We hope that our Theorem~\ref{thm:main-iso} can be used to show that their category of Koszul-perverse coherent sheaves provides a monoidal cluster categorification of the cluster algebras associated to quiver gauge theories.

\subsection*{Notations.}

For the duration of the paper we work with quantized K-theoretic Coulomb branch algebras $\Acr_{G,N}^q$ and quantized universally Laurent algebras $\Lbb_\quiver^q$. To lighten our notations we drop the super-script $q$ and write
$$
\Acr_{G,N} = \Acr_{G,N}^q, \qquad \Lbb_\quiver = \Lbb_\quiver^q.
$$

\subsection*{Acknowledgements.}

We are very grateful to Murad Alim, Mikhail Bershtein, Alexander Braverman, Michele Del Zotto, Philippe Di Francesco, Tudor Dimofte, Ilya Dumanski, Michael Finkelberg, Davide Gaiotto, Michael Gekhtman, Alexander Goncharov, Rinat Kedem, Daniil Klyuev, Vasily Krylov, Hiraku Nakajima, Michael Shapiro, Alexander Tsymbaliuk, Alex Weekes, and Harold Williams for many helpful discussions. G.S.\ has been supported by the NSF Standard Grant DMS-2302624. A.S.\ has been supported by the European Research Council under the European Union's Horizon 2020 research and innovation programme under grant agreement No 948885 and by the Royal Society University Research Fellowship.

\section{Coulomb branches of quiver gauge theories}
\label{sec:bfn}

\subsection{Quiver Coulomb branch algebras.}
\label{subsec:coulomb-def}

We start by recalling the construction of the quantized $K$-theoretic Coulomb branch algebra of a $3d$ $\Nc=4$ quiver gauge theory in the setup of Braverman, Finkelberg, and Nakajima, see~\cite{Nak16, BFN18}. Such a theory is determined by a \emph{weighted quiver}  $\Gamma$ with nodes $\Gamma_0$ and arrows $\Gamma_1$. The set of nodes
$$
\Gamma_0 = \Gamma_0^{\rg} \sqcup \Gamma_0^\rf.
$$
is partitioned into two subsets, which we refer to respectively as `gauge' and `flavor' nodes, with the property that {each flavor node is of valency 1, and there are no arrows connecting flavor nodes}.  Hence the set of arrows in $\Gamma$ decomposes as
$$
\Gamma_1 = \Gamma^\rg_1 \sqcup \Gamma^{\rf_+}_1 \sqcup \Gamma^{\rf_-}_1, 
$$
where $\Gamma^\rg_1$ is the set of all arrows $a$ such that both the source $s(a)$ and target $t(a)$ are gauge nodes, $ \Gamma^{\rf_+}_1$ the set of arrows whose source $s(a)$ is a flavor node, and $\Gamma^{\rf_-}_1$ those arrows for whose target $t(a)$ is a flavor node.
Each gauge node $i\in \Gamma_0^\rg$ is `weighted' with a finite dimensional complex vector space $V_i$, while at each flavor node $k\in\Gamma_0^\rf$ sits a 1-dimensional complex vector space $W_k$. We set
$$
V = \bigoplus_{i \in \Gamma^\rg_0} V_i \qquad\text{and}\qquad W = \bigoplus_{i \in \Gamma^{\rf}_0} W_i.
$$
and write
$$
\dimv_i = \dim(V_i), \quad i\in \Gamma^\rg_0.
$$
The gauge group
\begin{align}
\label{eq:gg}
G_V = \prod_{i \in \Gamma^\rg_0} GL(V_i)
\end{align}
is represented on the space
{
$$
N_\Gamma = \bigoplus_{a\in \Gamma_1} N_a,
$$
where
$$
N_a =
\begin{cases}
\Hom_\C(V_{s(a)}, V_{t(a)}) &\text{if} \quad a\in \Gamma_1^\rg, \\
\Hom_\C(W_{s(a)}, V_{t(a)}) &\text{if} \quad a\in \Gamma_1^{\rf_+}, \\
\Hom_\C(V_{s(a)}, W_{t(a)}) &\text{if} \quad a\in \Gamma_1^{\rf_-}.
\end{cases}
$$
}
We consider the torus
$$
T_W= \prod_{i \in \Gamma^\rf_0} GL(W_i) \simeq (\C^\times)^{\hm{\Gamma^\rf_0}}
$$
which is also naturally represented on the vector space $N_\Gamma$, and whose action commutes with that of the gauge group $G_V$ from~\eqref{eq:gg}.
We denote the defining characters for the factors in the torus $T_W$ by~$z_{k}$, $k\in \Gamma^\rf_0$, so that
$$
K_{T_W}(\pt) = \mathbb{Z} \big[z_{k}^{\pm1} \,\big|\, k \in \Gamma^\rf_0\big].
$$
There is also a further symmetry of $N_\Gamma$ coming from the torus
$$
T_H = H^1(\Gamma,\C^\times) \simeq (\C^\times)^g.
$$
In order to speak about this action concretely, we  fix a collection of 1-cochains $(\theta_1, \ldots, \theta_g)$ 
$$
\theta_i \colon \Z^{\Gamma_1} \longrightarrow \Z
$$
on $\Gamma$ whose cohomology classes form a basis for  $H^1(\Gamma,\mathbb{Z})$. Without loss of generality we may assume that $\theta_i(a)=0$ for any edge $a$ joining a gauge and a framing node. We then define the action of $(\mathbb{C}^\times)^g$ on $N_\Gamma$ by declaring that the $i$-th $\C^\times$ factor acts on the vector space $\Hom(V_j,V_k)$ associated to an arrow  $a \colon j \rightarrow k$ with the character $\theta_i(a)$. We denote the equivariant parameter corresponding to the class $[\theta_i]$ by $u_i$, so that
$$
K_{T_H}(\pt) \simeq \mathbb{Z}\hs{u_1^{\pm 1}, \ldots, u^{\pm1}_g},
$$
and write 
$$
u(a) = \prod_{i=1}^g u_i^{{-}\theta_i(a)}
$$ 
for any edge $a\colon j \rightarrow k$. Note that if $a$ is a loop connecting a gauge node to itself we must necessarily have $u(a)\neq1$ by virtue of the $\theta_i$ forming a cohomology basis. Finally, we define the \emph{full torus of flavor symmetry} to be
$$
T_F = T_W \times T_H \simeq \Z\big[z_k^{\pm1},u_a^{\pm1} ~\big|~ k \in \Gamma_0^\rf, \, 1 \le a \le g\big].
$$

Now let $G,G_F$ be a pair of complex reductive groups and $N$ a complex representation of $G\times G_F$.
We consider the balanced product
\begin{align}
\label{eq:TGN}
\mathscr{T}_{G,N} = G\brr{z} \times_{G\bss{z}} N\bss{z} = (G\brr{z} \times N\bss{z})/G\bss{z}
\end{align}
where the quotient is taken with respect to the action $h \circ (g,n) = (gh^{-1}, hn)$. The space~\eqref{eq:TGN} can be thought of as an infinite-rank $(G\times G_F)\bss{z}$-equiariant vector bundle over the affine Grassmannian 
$$
\Gr_{G}=G\brr{z}/{G\bss{z}}
$$
via the projection to the first factor in the balanced product. Following~\cite{BFN18}, we consider the \emph{variety of triples}
$$
\Rcr_{G,N} = \hc{[g,s] \in G\brr{z} \times_{G\bss{z}} N\bss{z} \,\big|\, gs \in N\bss{z}}. 
$$
The space $\Rcr_{G,N}$ is preserved by the action of the arc group $(G\times G_F)\bss{z}$, and carries an action of $\C^\times$ by loop rotation, that is by rescaling of the loop coordinate $z$. Again following~\cite{BFN18, FT17} we consider the double cover $\widetilde\C^\times$ of $\C^\times$, writing $q^2$ for the equivariant parameter corresponding to $\C^\times$, and $q$ for the one corresponding to $\widetilde{\C}^\times$, so that
$$
K_{\C^\times}(\pt) = \Z \hs{q^{\pm2}} \qquad\text{and}\qquad K_{\widetilde\C^\times}(\pt) = \Z \hs{q^{\pm1}}.
$$
Following Section 2(ii) and Remark 3.9 in~ \cite{BFN18}, we recall the following definition.

\begin{defn}
The \emph{quantized $K$-theoretic Coulomb branch ring} associated to the data $(G,G_F,N)$ is defined as the equivariant algebraic $K$-theory
$$
\Acr_{G,G_F,N} = K_{(G\times G_F)\bss{z} \rtimes \widetilde{\C}^\times} \hr{\Rcr_{G,N}}
$$
of the variety of triples $\Rcr_{G,N}$.
\end{defn}

We recall from~\cite[Section 3(vii)]{BFN18} that the space $\Rcr_{G,N}$ fits into a convolution diagram similar to that for the affine Grassmannian, and that this convolution endows the abelian group $\Acr_{G,G_F,N}$ with the structure of an associative algebra over the commutative ring $K_{G_V \times T_F \times \widetilde{\C}^\times}(\pt)$. We also recall that
$$
\Acr_{G,G_F,0} \simeq \Acr_{G,0} \otimes_\Z K_{G_F}(\pt),
$$
where we abbreviate $\Acr_{G,0} = \Acr_{G,\{e\},0}$.

As in~\cite{BFN19}, we can apply this construction to the data $(G_V, T_F, N_\Gamma)$ associated to a weighted quiver $\Gamma$ to produce a variety of triples $\Rcr_{\Gamma}=\Rcr_{G_V,N_\Gamma}$ and a convolution algebra
$$
\Acr_{\Gamma} = K_{(G_V\times T_F)\bss{z} \rtimes \widetilde{\C}^\times } \hr{\Rcr_{\Gamma}}
$$
which we refer to as the \emph{$K$-theoretic Coulomb branch ring}  of the quiver theory.

\begin{remark}
\label{rem:G-to-T-sym}
In a variant of the quiver gauge theory setup described above, one can drop the requirement that the spaces $W_k$ at the framing nodes of quiver $\Gamma$ be 1-dimensional and consider the flavor group
$$
G_F = G_W \times T_H, \qquad G_W = \prod_{i \in \Gamma^\rf_0} GL(W_i).
$$
This seemingly more general picture is related to the one described here as follows: let $T_W \subset G_W$ be the maximal torus, and denote by
$$
S_W = \prod_{i \in \Gamma^\rf_0} S_{\dim(W_i)}
$$
the Weyl group of $G_W$. Then we have
$$
\Acr_{G_V,G_F,N_\Gamma} \simeq \Acr_{G_V,T_F,N_\Gamma}^{S_W} \simeq \Acr_{\Gamma'}^{S_W}
$$
where the quiver $\Gamma'$ is obtained from $\Gamma$ via replacing each node carrying the space $W_k$ with $k$ nodes carrying the 1-dimensional weight subspaces of $W_k$.
\end{remark}

\subsection{Minuscule monopole operators}
\label{monopole-sect}

Recall that for a reductive group $G$, the affine Grassmannian
$$
\Gr_G = G\brr{z}/G\bss{z}
$$
is stratified by finite-dimensional $G\bss{z}$-orbits:
$$
\Gr_G = \bigsqcup_{\lambda \in P^+} \Gr_G^\la, \qquad \Gr_G^\la = G\bss{z}z^\lambda,
$$
where $P^+$ is the dominant cone in the coweight lattice of $G$, or equivalently the weight lattice for the dual group $G^\vee$. An orbit $\Gr_G^\la$ is closed if and only if the dominant coweight $\la$ is minuscule, i.e.\ all weights in the corresponding irreducible representation $L_\lambda$ of the dual group $G^\vee$ are conjugate under the Weyl group.

Given a dominant minuscule coweight $\la$ of $G_V$, denote by $\Rcr_\la$ the pre-image of the corresponding closed orbit under the canonical projection from $\Rcr_{\Gamma}$ to $\Gr_{G_V}$. Then the $K$-theory class $\hs{\Oc_{\Rcr_\la}}$ of its structure sheaf defines an element of $\Acr_\Gamma$ called a \emph{minuscule monopole operator.}

{For notational convenience, we opt to identify both the weight and the coweight lattices of $GL(V_i)\simeq GL_{\dimv_i}$ with the lattice $\mathbb{Z}^{\dimv_i}$, and write $\{\eps_{i,r}\}_{r=1}^{\dimv_i}$ for the standard basis in the latter, with respect to which the invariant bilinear form reads $\ha{\eps_{i,r},\eps_{j,s}} = \delta_{i,r}\delta_{j,s}$.} We denote by
$$
\varpi_{i,n} = \sum_{r=1}^n\eps_{i,r} 
$$
the $n$-th fundamental coweight of $GL(V_i)$, and set 
$$
\varpi_{i,n}^* = -w_0(\varpi_{i,n}) = \varpi_{i,\dimv_i-n}-\varpi_{i,\dimv_i}. 
$$

The closed orbit  $\Gr_{G_V}^{\varpi_{i,n}}$ is isomorphic to the finite dimensional Grassmannian $\Gr(V_i,n)$ of $n$-dimensional quotients of $V_i$. We write $\Qc_{\varpi_{i,n}}$ for the pullback to $\Rcr_{\varpi_{i,n}}$ of the tautological rank $n$ vector bundle of quotients on $\mathrm{Gr}(V_i,n)$, and $\Sc_{{\varpi_{i,n}}}$ for the tautological rank $\dimv_i-n$ bundle of kernels. For a vector bundle $\mathcal{V}$, we set $e_p(\mathcal{V}) = \bigwedge^p(\mathcal{V})$, so that we can evaluate arbitrary Laurent polynomials
$$
f \in K_{GL_{\dimv_i}}(\Gr(V_i,n)) \simeq K_{T_{\dimv_i}}(\pt)^{S_n \times S_{\dimv_i-n}},
$$
symmetric under the stabilizer of ${\varpi_{i,n}}$, on the bundles $\Qc_{\varpi_{i,n}}$, $\Sc_{{\varpi_{i,n}}}$, where we write $T_{\dimv_i}$ for the maximal torus of $GL_{\dimv_i}$, producing $K$-theory classes 
$$
f(\Qc_{{\varpi_{i,n}}},\Sc_{{\varpi_{i,n}}})\in K_{(G_V\times T_F)\bss{z} \rtimes \widetilde{\C}^\times }(\Rcr_{\varpi_{i,n}}).
$$
Similarly, the closed orbit $\Gr_{G_V}^{\varpi^*_{i,n}}$ is isomorphic to the finite dimensional Grassmannian $\Gr(n,V_i)$ of $n$-dimensional subspaces of $V_i$. We denote by $\Sc_{\varpi^*_{i,n}},\Qc_{\varpi^*_{i,n}}$ the pull-back to $\Rcr_{\varpi_{i,n}^*}$ of the tautological sub- and quotient bundles on $\Gr_{G_V}^{\varpi^*_{i,n}}$, with equivariant structure twisted by tensoring with the defining character $q^2$ of the loop rotation torus $\mathbb{C}^\times$.

Given a collection of integers
\beq
\label{eq:n-col}
\bs n = \hc{n_i \,|\, 0 \le n_i \le \dimv_i}_{i \in \Gamma^\rg_0},
\eeq
we define a pair of dominant minuscule coweights
\beq
\label{min-cow}
\varpi_{\bs n} = \sum_{i \in \Gamma^\rg_0} \varpi_{i,n_i} \qquad\text{and}\qquad \varpi_{\bs n}^* = \sum_{i \in \Gamma^\rg_0} \varpi^*_{i,n_i}.
\eeq
Then further fixing an element
$$
\bs f  = \hc{f_i \,|\, i \in \Gamma^\rg_0} \in \bigotimes_{i\in\Gamma_0^\rg}K_{T_{\dimv_i}}(\pt)^{S_{n_i}\times S_{\dimv_i-n_i}}
$$
we consider the \emph{dressed minuscule monopole operators}
$$
E_{\bs n, \bs f} = \hs{\bigotimes_{i \in \Gamma^\rg_0} f_i\Big(\Qc_{\varpi_{i,n_i}}, \Sc_{\varpi_{i,n_i}}\Big)}, \qquad 
F_{\bs n, \bs f} = \hs{\bigotimes_{i \in \Gamma^\rg_0} f_i\Big(\Sc_{\varpi^*_{i,n_i}},\Qc_{\varpi^*_{i,n_i}}\Big)}.
$$
In particular, we have
$$
E_{\bs n, \bs 1} = \hs{\Oc_{\Rcr_{\varpi_{\bs n}}}} \qquad\text{and}\qquad F_{\bs n, \bs 1} = \hs{\Oc_{\Rcr_{\varpi^*_{\bs n}}}},
$$
where $\bs 1_i = 1$ for all $i \in \Gamma_0^\rg$.

The following result was proved for cohomological Coulomb branches in~\cite[Proposition 3.1]{Wee19}, following the general procedure outlined in~\cite[Proposition 6.8]{BFN18}. The proof in the $K$-theory case we consider here goes the same way, and is written in~\cite[Proposition 4.7]{DK25}.

\begin{prop}
\label{prop-gen}
As an algebra over the ring $K_{G_V \times T_F \times \widetilde{\C}^\times}(\pt)$,
the quantized Coulomb branch $\Acr_\Gamma$ is generated by the set of all dressed minuscule monopoles  
$
E_{\bs n, \bs f}$ and $F_{\bs n, \bs f}$ as $\bs n$ ranges over all possible collections of the form~\eqref{eq:n-col}.
\end{prop}

Since we will refer to it later, we briefly recall following~\cite{Wee19} the combinatorics behind Proposition~\ref{prop-gen}. The set of all weights for the maximal torus
$$
T_V = \otimes_{i \in \Gamma_0} T_{\dimv_i} \subset G_V
$$
appearing in the representation $N_\Gamma$ gives an arrangement of hyperplanes in the Lie algebra $\tgt_V = \mathrm{Lie}(T_V)$, which further divide the dominant Weyl chamber in $\tgt_V$ into a union of finitely many chambers~$\mathfrak{C}$. The semigroup of integral points of each chamber $\mathfrak{C}$ is generated by some finite collection of dominant coweights, and the corresponding set of all dressed minuscule monopoles as we range over all chambers gives a set of generators for $\Acr_\Gamma$, see \cite{BFN18}. In fact, we can obtain a (possibly redundant) set of generators by further refining the hyperplane arrangement so that it contains all the hyperplanes $\{\eps_{i,m}=\eps_{j,n}\}$ and $\{\eps_{i,m}=0\}$, where we have identified the weight and coweight lattices. The chambers $\mathfrak C$ of the refined arrangement are obtained by intersecting Weyl chambers of the group $GL(V)$ with hyperplanes $\eps_{i,n}=0$, and take the form
\beq
\label{eq:ref-Weyl-chamber}
\mathfrak C = \hc{\eps_{i_1,r_1} \ge \ldots \ge \eps_{i_k,r_k} \ge 0 \ge \eps_{i_{k+1},r_{k+1}} \ge \ldots \ge \eps_{i_{\dim(V)},r_{\dim(V)}}},
\eeq
where $\dim(V) = \sum_{i \in \Gamma_0^\rg} \dimv_i$. Note that lattice of integral points in the refined Weyl chambers $\mathfrak C$ subdividing the dominant Weyl chamber of $GL(V)$ are generated by the dominant minuscule coweights~\eqref{min-cow}, so the claim of the Proposition follows.

\subsection{Equivariant localization}
\label{subsec:localization}

Using Proposition~\ref{prop-gen} and the localization formula in equivariant $K$-theory it is possible to give an explicit description of $\Acr_\Gamma$ as a subalgebra inside a ring of $q$-difference operators.  Let $w_{i,r}$ be the character of the torus $T_V$ corresponding to the weight $\eps_{i,r} \in \tgt_V^*$, so that
$$
K_{T_V \times T_F \times \widetilde{\C}^\times}(\pt) \simeq \Z\hs{q^{\pm1}, w_{i,r}^{\pm1}, z_{k}^{\pm1}, u_a^{\pm1}},
$$
where $i \in \Gamma^\rg_0$, $1 \le r \le \dimv_i$, $k\in \Gamma^\rf_0$, and $1 \le a \le g$. Then we have
$$
K_{G_V\times T_F \times \widetilde{\C}^\times}(\pt) \simeq {K_{T_V \times T_F \times \widetilde{\C}^\times}(\pt)}^{S_V},
$$
where
$$
S_V = {\prod_{i \in \Gamma^\rg_0}S_{\dimv_i}},
$$
and the symmetric group $S_{\dimv_i}$ acts by permuting the variable set $\bs w_{i,\bullet} = \hc{w_{i,r} \,|\, 1 \le r \le \dimv_i}$.

Now we consider the Coulomb branch ring in which by we replace the group $G_V$ by its maximal torus $T_V$, and the representation $N_\Gamma$ by the zero representation. As was shown in~\cite{BFN18}, the resulting algebra $\Acr_{T_V,T_F,0} = K_{(T_V \times T_F)\bss{z} \rtimes \widetilde{\C}^\times }\hr{\Rcr_{T_V,0}}$ is isomorphic to the algebra
$$
\Dc(\Gamma) = \Dc(T_V) \otimes_{\mathbb{Z}} K_{T_F}(\pt),
$$
where
$$
\Dc(T_V) = \Z[q^{\pm1}]\ha{D_\mu,w_\la \,|\, \mu\in\Prm, \la \in \Prm^\vee},
$$
is the ring of $q$-difference operators on the torus $T_V$, subject to the defining relations
$$
D_\mu w_\la = q^{2\ha{\mu,\la}}w_\la D_\mu,
$$
and we write $\Prm$, $\Prm^\vee$ for the weight and coweight lattices of $G_V$ respectively. We also set $D_{i,r} = D_{\eps_{i,r}}$, $w_{j,s} = w_{\eps_{j,s}}$, so that for $\mu = \sum_{i,r} \mu_{i,r}\eps_{i,r}$ we have
$$
D_\mu = \prod_{i\in \Gamma_0^\rg}\prod_{r=1}^{\dimv_i}D_{i,r}^{\mu_{i,r}},
\qquad
w_\mu = \prod_{i\in \Gamma_0^\rg}\prod_{r=1}^{\dimv_i}w_{i,r}^{\mu_{i,r}},
$$
and
$$
\Dc(T_V) \simeq \Z[q^{\pm1}]\ha{D_{i,r},w_{i,r}}/\langle D_{i,r}w_{j,s} = q^{2\delta_{i,j}\delta_{r,s}} w_{j,s}D_{i,r}\rangle.
$$
We also write 
$$
\Delta_+(G_V) = \hc{\eps_{i,r}-\eps_{i,s}~\big|~ i\in \Gamma_0^\rg, \, 1\leq r<s\leq \dimv_i}, \qquad \Delta_-(G_V) = -\Delta_+(G_V)
$$
for the set of all positive roots of the gauge group $G_V$, and
$$
\Delta(G_V) =\Delta_+(G_V) \sqcup \Delta_-(G_V) 
$$
for the set of all roots. Hence a root $\alpha = \eps_{i,r}-\eps_{i,s}\in\Delta(G_V)$ corresponds to the class
$$
w_\alpha = w_{i,r}w^{-1}_{i,s} \in K_{T_V}(\pt).
$$

The localization theorem in equivariant $K$-theory (followed by the restriction to the zero section in the vector bundle~\eqref{eq:TGN} for the abelian theory, see Appendix A of \cite{BFN19} or Remark 5.23 of~\cite{BFN18})  provides an injective algebra homomorphism
\begin{align}
\label{eq:localization-hom}
\mathbf{z}^*\hr{\iota_*}^{-1} \colon \Acr_\Gamma \hookrightarrow \Dc_\loc(\Gamma)
\end{align}
of the quantized Coulomb branch $\Acr_\Gamma$ into the localization $\Dc_\loc(\Gamma)$ of $\Dc(\Gamma)$ at the Ore denominator set
$$
Ø(\Gamma) = \hc{1 -q^{2k}w_{\alpha} ~\Big|~ \alpha\in \Delta_+(G_V), ~k\in\mathbb{Z}}.
$$
Since the corresponding orbit closures are smooth with isolated $T_V$-fixed points, images of the dressed minuscule monopole operators under the embedding~\eqref{eq:localization-hom} can be easily computed using Koszul resolutions for the structure sheaves of the fixed points, see~\cite{BFN19, FT17}. In particular, the images of the generators $E_{\bs n, \bs f}, F_{\bs n, \bs f}$ of $\Acr_\Gamma$ from Proposition~\ref{prop-gen} are given by

\beq
\label{Ebs}
\begin{aligned}
\mathbf{z}^*\hr{\iota_*}^{-1}(E_{\bs n, \bs f}) = \sum_{\bs J:|J_l|=n_l} \prod_{\substack{a:i\to j\\ a\in \Gamma_1^\rg}} \prod_{\substack{r \in J_i \\ s \notin J_j \\ (j,s) \ne (i,r)}} \hr{1-qu(a)w_{i,r}w_{j,s}^{-1}} 
\prod_{\substack{a:i\to k\\ a\in \Gamma_1^{\rf_-}}}\prod_{r\in J_i}\hr{1-qw_{i,r}z_{k}^{-1}}\\
\prod_{l \in \Gamma^\rg_0} f_l(\bs w) \prod_{\substack{r \in J_l \\ s \notin J_l}} \frac{1}{1-w_{l,s}w_{l,r}^{-1}} \prod_{r \in J_l} D_{l,r},
\end{aligned}
\eeq

\beq
\label{Fbs}
\begin{aligned}
\mathbf{z}^*\hr{\iota_*}^{-1}(F_{\bs n, \bs f}) = \sum_{\bs J:|J_l|=n_l} \prod_{\substack{a:j\to i\\ a\in \Gamma_1^\rg}} \prod_{\substack{r \in J_i \\ s \notin J_j \\ (j,s) \ne (i,r)}} \hr{1-qu(a)w_{j,s}w_{i,r}^{-1}}
\prod_{\substack{a:k\to j\\ a\in \Gamma_1^{\rf_+}}}\prod_{r\in J_j}\hr{1-qz_{k}w_{j,r}^{-1}}\\
\prod_{l \in \Gamma^\rg_0} f_l(\bs w) \prod_{\substack{r \in J_l \\ s \notin J_l}} \frac{1}{1-w_{l,r}w_{l,s}^{-1}} \prod_{r \in J_l} D_{l,r}^{-1}.
\end{aligned}
\eeq
The summation in the formulas above is taken over all tuples of subsets $\bs J = (J_l \,|\, l \in \Gamma^\rg_0)$ indexed by the gauge nodes of $\Gamma$, where $J_l \subset \hc{1, \dots, d_l}$ and $|J_l|=n_l$. We also  write $r \notin J_l$ for $r \in \hc{1, \dots, d_l} \setminus J_l$. Note that the condition $(j,s) \ne (i,r)$ is void unless the arrow $a$ is a 1-cycle. 

We can rewrite these formulas in a slightly less cumbersome form as follows. Denote by
$$
\supp(N_\Gamma) \subset \mathfrak t_{V}^* \otimes \mathfrak t_{F}^*
$$
the set of all weights for the torus $T_V\times T_F$ appearing in the representation $N_\Gamma$. Given a weight $\chi \in \supp(N_\Gamma)$, we write $\chi = (\chi_V, \chi_W, \chi_H)$ for its components with respect to the tori $T_V$, $T_W$, and $T_H$. Note that the multiplicities of these weights are all 1, since we work with the full torus of flavor symmetry. Then for $\la$ of the form~\eqref{min-cow}, we have
\beq
\label{eq:zfO}
\mathbf{z}^*\hr{\iota_*}^{-1}(f\otimes[\Oc_{\Rcr_{\la}}]) = \sum_{[\sigma]\in S_V/\St_\la}\sigma\hr{ f(\bs w)\prod_{\substack{\chi\in\supp(N_\Gamma)\\ \ha{\chi,\la}<0}}(1-qX_\chi^{-1})\prod_{\substack{\alpha\in \Delta(G_V) \\ \ha{\alpha,\la}<0 }}\frac{1}{1-w_\alpha}D_{\la}}
\eeq
where we write
$$
f \otimes \hs{\Oc_{\Rcr_\la}} =
\begin{cases}
E_{\bs n, \bs f}, &\text{if}~\la = \varpi_{\bs n}, \\
F_{\bs n, \bs f}, &\text{if}~\la = \varpi^*_{\bs n},
\end{cases}
\qquad\qquad
X_\chi = w_{\chi_V}z_{\chi_W}u_{\chi_H},
$$
and the sum is taken over coset representatives for the stabilizer $\St_\la$ of the coweight $\la$ in the product of symmetric groups $S_V$.

We can further rewrite this latter formula in terms of conjugation by a formal power series
 \begin{align}
\label{eq:qdl-def}
\Psi(X)= \prod_{n=0}^\infty (1-q^{2n+1}X)^{-1} \in \mathbb{Q}(q)[[X]]\subset \mathbb{Z}\brr{q}[[X]]
\end{align}
known as the \emph{(compact) quantum dilogarithm function}, which satisfies the $q$-difference equation
\begin{align}
\label{eq:qde}
\Psi(q^2X) = (1-qX)\Psi(X).
\end{align}
Given a dominant coweight $\la$, we define formal series
\begin{align}
\label{eq:gauge-qdl}
\mathbf{\Psi}_{G_V;\la} &= \prod_{\substack{\alpha\in \Delta(G_V) \\ \ha{\alpha,\la}\le0}}\Psi(qw_\alpha), \\
\mathbf{\Psi}_{N_\Gamma;\la} &= \prod_{\substack{\chi\in \supp(N_\Gamma) \\ \ha{\chi,\la}\le0 }}\Psi(X_\chi^{-1}).
\end{align}

For dominant $\la$, the series $\mathbf{\Psi}_{G_V;\la}$ takes values in the completion of $K_{T_V\times T_F}(\pt)\otimes_\Z \mathbb{Q}(q)$ associated to the cone $\Delta_-(G_V)$ of negative roots, while $\mathbf{\Psi}_{N_\Gamma;\la}$ lives in the completion corresponding to the cone in the character lattice dual to the dominant cone in the cocharacter lattice. Then using the $q$-difference equation~\eqref{eq:qde} we can write~\eqref{eq:zfO} in the compressed form
\begin{align}
\label{eq:compressed}
\mathbf{z}^*\hr{\iota_*}^{-1}(f\otimes[\Oc_{\Rcr_{\la}}]) = \sum_{[\sigma]\in S_V/\St_\la}\sigma \hr{\Ad_{\mathbf{\Psi}^{-1}_{N_\Gamma;\la}}\circ \Ad_{\mathbf{\Psi}^{-1}_{G_V;\la}}(f(\bs w)D_{\la})}.
\end{align}

\begin{remark}
Let $L$ be a 1-dimensional representation of the torus $T=\C^\times$ and write $X\in K_T(\pt)$ for its class in the representation ring of $T$. Then the vector space $L[z]$ of all regular maps from $\mathbb{A}^1$ to $L$ is an infinite dimensional representation of the torus $T\times \widetilde{\mathbb{C}}^\times$, where following the convention of~\cite{BFN18} we define the action by $(t,s)\cdot \phi = t\cdot(s\phi(s^2z))$. In the completed representation ring $\widehat{K}_{T\times \widetilde{\mathbb{C}}^\times}(\pt)=\mathbb{Z}\brr{q}[X^{\pm1}]$, the class of this representation is
$$
L[z] = \frac{qX}{1-q^2} \in \widehat{K}_{T\times \widetilde{\mathbb{C}}^\times}(\pt),
$$
and so the quantum dilogarithm series~\eqref{eq:qdl-def} represents the class of the symmetric algebra of the representation $L[z]$:
$$
\Psi(X) = S^\bullet\left(L[z] \right) \in \widehat{K}_{T\times \widetilde{\mathbb{C}}^\times}(\pt).
$$
\end{remark}

\subsection{Coulomb branch algebra of the pure $GL(V)$ gauge theory}
\label{subsec:residue-sec}
 In this section we describe the Coulomb branch algebra 
$$
\Acr_{G_V,0} = \Acr_{G_V,\hc{e},0} = K_{G_V \bss{z} \rtimes \widetilde{\C}^\times }\hr{\Gr_{G_V}},
$$
associated to gauge group $G_V=\prod_i GL_{\dimv_i}$, trivial flavor group, and the zero representation. As in Section~3(vii) of~\cite{BFN18}, we have
$$
\Acr_{G_V,0} \simeq \bigotimes_{i\in \Gamma_0^\rg}\Acr_{\mathrm{GL}_{\dimv_i},0},
$$
where the tensor products are taken over $\mathbb{Z}[q^{\pm1}]$ and the tensor factors commute. So it suffices to understand the $K$-theoretic Coulomb branch algebra $\Acr_{\mathrm{GL}_{\dimv},0}$ associated to the pure $GL_\dimv$ gauge theory, which following~\cite{SS25} we do by characterizing its image under the localization embedding~\eqref{eq:localization-hom} in the localized ring of $q$-difference operators
$$
\label{eq:Dc-loc}
\Dc_\loc(T_\dimv) = \Z[q^{\pm1}]\ha{D^{\pm1}_{r},w^{\pm1}_{r}}_{r=1}^{\dimv}\big[(1-q^{2k}w_\alpha)^{-1}~\big|~\alpha\in\Delta_+(GL_\dimv), \, k\in\mathbb{Z}\big].
$$
In this case, the formulas~\eqref{Ebs} and~\eqref{Fbs} for the images of the dressed minuscule monopoles simplify:
\begin{align}
\label{Eims}
\mathbf{z}^*\hr{\iota_*}^{-1}(E_{n,f}) &= \sum_{|J|=n} f(\bs w) \prod_{\substack{r\in J \\ s\notin J}} \frac{1}{1-w_{s}w_{r}^{-1}} \prod_{j \in J} D_j, \\
\label{Fims}
\mathbf{z}^*\hr{\iota_*}^{-1}(F_{n,f}) &= \sum_{|J|=n} f(\bs w) \prod_{\substack{r\in J \\ s\notin J}} \frac{1}{1-w_{r}w_{s}^{-1}} \prod_{j \in J} D_j^{-1}.
\end{align}

In particular, let us write
$$
\bs w_J = \hc{w_j \,|\, j \in J}, \qquad \bs w_{J^c} = \hc{w_j \,|\, j \notin J}.
$$
Then for $f(\bs w) = f_1(\bs w_J) f_2(\bs w_{J^c})$ and $f^\star(\bs w) = f_2(\bs w_J) f_1(\bs w_{J^c})$ we have
$$
F_{n,f} \cdot E_{\dimv,1} = E_{\dimv-n,f^{\star}}, \qquad E_{n,f} \cdot F_{\dimv,1} = F_{\dimv-n,f^\star}.
$$

The algebra $\Dc_\loc(T_\dimv)$ carries an internal grading by the coweight lattice $\Z^{\dimv}$, whose degree $\mu=(\mu_1,\ldots, \mu_\dimv)$ piece is given by
$$
A\in \Dc_\loc^\mu(T_\dimv) \iff Aw_i = q^{2\mu_i}w_iA, \quad 1\leq i\leq \dimv.
$$ 
In concrete terms, $\Dc_\loc^\mu(T_\dimv)$ consists of all elements of the form $h(\bs w,q) D_\mu$, where $h(\bs w,q)$ is a rational function in variables $\bs w$ and $q$ regular outside the union of the divisors
$$
d_{\alpha,k} = \hc{q^{2k}w_\alpha=1}.
$$
Given a general element $A\in \Dc_\loc(T_\dimv)$ we write $A_\mu$ for its degree $\mu$ component, so that
$$
A = \sum_{\mu\in\mathbb{Z}^\dimv} A_\mu, \quad A_\mu = h_\mu(\bs w,q)D_\mu \in \Dc_\loc^\mu(T_\dimv) .
$$

Recall that the symmetric group ${S_\dimv}$ acts on $\Dc_\loc(T_\dimv)$ via permutation of the canonically conjugate pairs $(D_r,w_r)$. We write $\Dc_\loc(T_\dimv)^{S_\dimv}$ for the subalgebra of ${S_\dimv}$-invariant elements in $\Dc_\loc(T_\dimv)$. We denote by $s_{\alpha}$ the reflection associated to a root $\alpha\in \Delta_+$, where $\Delta_+ = \Delta_+(GL_\dimv)$ for the remainder of this section. The following definition was introduced in~\cite{SS25}. It is modelled on the constructions of~\cite{GKV97}, in particular their Definition~1.3 and Theorem~1.4.

\begin{defn}
\label{def:Dres}
We define $\Dres(T_\dimv)$ to be the subspace of $\Dc_\loc(T_\dimv)^{S_\dimv}$ consisting of elements 
$$
A =\sum_{\mu\in\Z^\dimv} h_\mu(\bs w,q)D_\mu, \quad 
$$ such that 
\begin{enumerate}
\item for each $\mu \in \mathbb{Z}^\dimv$, $\alpha\in \Delta_+$, and $k \in \Z$ the rational function $h_\mu(\bs w,q)$ has at worst a simple pole at $d_{\alpha,k}$;
\item for all $\mu \in \mathbb{Z}^\dimv$, $\alpha\in\Delta_+$, and $k \in \Z$ the corresponding residues satisfy
\beq
\label{eq:residue-condition}
\res_{\alpha,k}\hr{h_{\mu}+h_{s_\alpha(\mu)+k\alpha}}=0.
\eeq
\end{enumerate}
\end{defn}

\begin{remark}
\label{rmk:minuscule-res}
If an element $A=\sum_\mu A_\mu\in  \Dc_\loc(T_\dimv)^{S_\dimv}$ has at worst simple poles and lies in the direct sum of graded pieces $\Dc^\mu_\loc(T_\dimv)$ where $\mu$ runs over all (not necessarily dominant) minuscule coweights, then it is easy to see that the residue condition~\eqref{eq:residue-condition} is automatically satisfied. 
\end{remark}
\begin{prop}[Lemmas 6.6, 6.12, and 6.14(i) in~\cite{SS25}] \,
\label{prop:ss25}
\begin{enumerate}
\item $\Dres(T_\dimv)$ forms a subalgebra in $\Dc_\loc(T_\dimv)$;
\item The subalgebra $\Dres(T_\dimv)$ is generated by its intersection with the direct sum of graded pieces $\Dc^\mu_\loc(T_\dimv)$ where $\mu$ runs over all minuscule coweights (not necessarily dominant, and including 0).
\item The image of $\Acr_{GL_\dimv,0}$ in $\Dc_\loc(T_\dimv)$ under the localization embedding $(\iota_*)^{-1}$ coincides with $\Dres(T_\dimv)$.
\end{enumerate}
\end{prop}

\begin{proof}
The proof of part (1) is a straightforward calculation, see Lemma 6.6 of~\cite{SS25}.
Because we will use a similar argument later to prove Proposition~\ref{prop:conj-crit}, let us recall the proof of part (2) of the Proposition (which easily implies part 3). 
For a dominant coweight $\mu \in \mathbb{Z}^{\dimv}$, we define $\widetilde{D}_\mu \in \Dc_\loc(T_\dimv)$ to be the result of conjugating $D_\mu$ by the  element 
$\mathbf{\Psi}^{-1}_{GL_\dimv;\mu}$ from~\eqref{eq:gauge-qdl}:
$$
\widetilde{D}_\mu = \Ad_{\mathbf{\Psi}^{-1}_{GL_\dimv}}(D_\mu) = \frac{1}{d_\mu(\bs w)} D_\mu,
$$
where 
\begin{align}
\label{eq:denom}
 d_\mu(\bs w) =\prod_{\alpha \in \Delta_+} \prod_{k=1}^{\ha{\alpha,\mu}} \hr{1-q^{2(k-1)}w^{-1}_\alpha}.
\end{align}
Given a general element $A$ of $\Dres(T_\dimv)$, let us consider the finite set of summands $A_\mu$ of $A$ having maximal (in dominance order) internal degree as elements of $\Dc_\loc(T_\dimv)$. Since $A$ is Weyl-invariant, this gives us a finite set of dominant coweights $\mu$. Define the total \emph{height} of $A$ to be the sum of the heights $\ha{\rho,\mu}$ of each of these coweights, where $\rho$ is the Weyl vector, which is equal to the half sum of all positive roots. Now fix one such coweight $\mu$. The residue conditions~\eqref{eq:residue-condition} imply that the summand $A_\mu$ is regular except possibly at the divisors $d_{\alpha,k}$ with $0 \le k < \ha{\alpha,\mu}$, at which $A$ may have simple poles. Indeed, let us first assume that $k\ge\ha{\alpha,\mu}$, which implies $s_\alpha(\mu)+k\alpha \ge \mu$. Then the condition~\eqref{eq:residue-condition} together with the dominance-maximality of $A_\mu$ among the summands of $A$ implies that $\res_{\alpha,k}(h_\mu)=0$. On the other hand, if $k<0$ then $s_\alpha(\mu)+k\alpha=s_\alpha(\mu-k\alpha)$ belongs to the Weyl-orbit of the element $\mu-k\alpha > \mu$. By the Weyl-invariance of $A$ we have $A_{s_\alpha(\mu)+k\alpha}=s_\alpha(A_{\mu-k\alpha})=0$, and hence once again $\res_{\alpha,k}(h_\mu)=0$ by~\eqref{eq:residue-condition}. This shows the summand $A_\mu$ can indeed only have (simple) poles at the divisors $d_{\alpha,k}$ with $0 \le k<\ha{\alpha,\mu}$.

It follows that $A_\mu$ has the form $h(\bs w) \widetilde D_\mu$ for some $h \in \Z[q^{\pm1}][w_1^{\pm1}, \ldots, w_\dimv^{\pm1}]^{\St_\mu}$. To finish the proof, we will show that there exists a sum of products of dressed minuscule monopole operators whose leading term in dominance order coincides with $A_\mu$, and with all other summands having degree strictly less than $\mu$ in dominance order. Subtracting this sum of products of minuscule monopoles from $A$ we will obtain a new element $A'\in \Dres(T_\dimv)$ with strictly smaller total height, and so by induction it will follow that $\Dres(T_\dimv)$ is generated by the minuscule elements.

To this end, write $\mu'=(\mu'_1,\ldots, \mu'_l)$ for the partition conjugate to $\mu-\mu_d\varpi_d$, so that parts of $\mu'$ encode the expansion of $\mu$ into the basis of dominant weights:
$$
\mu = \mu_d\varpi_d + \sum_{i=1}^l \varpi_{\mu'_i}.
$$
Given a subset $J \subset \hc{1, \ldots, d}$, we write $S_J$ for the symmetric group $S_{\hm{J}}$ permuting the variable set $\hc{w_j \,|\, j \in J}$. Then the stabilizer of $\mu$ is given by
$$
\St_\mu \simeq \prod_{i=0}^{l}S_{\hc{\mu'_{i+1}+1, \ldots, \mu'_i}},
$$
where we understand $\mu'_0=d$ and $\mu'_{l+1}=0$. We need to show that any Laurent polynomial in $w_1,\ldots, w_d$ symmetric under $\St_\mu$ is expressible as a sum of products of the form $f_0\cdots f_l$, where $f_i = f_i(w_1,\ldots, w_{\mu'_i})$ is symmetric under the action of $S_{\mu_i'}$ permuting its arguments. In other words, we claim that the ring $K_G(\pt)^{\St_\mu}$ of $\St_\mu$-symmetric Laurent polynomials is generated by the subrings of $S_{\hc{1,\ldots,\mu_i'}}$-symmetric Laurent polynomials. To prove this, it suffices to show that each elementary symmetric polynomial
$$
e_{k}(w_{\mu'_{i+1}+1}, \ldots, w_{\mu'_i})
$$ of $K_G(\pt)^{\St_\mu}$ can be expressed in terms of elementary symmetric polynomials of the form $e_j(w_1,\ldots, w_{\mu'_r})$. But this follows by induction over $k$, since
$$
e_{k}(w_{\mu'_{i+1}+1}, \ldots, w_{\mu'_i}) = e_{k}(w_1,\ldots,w_{\mu'_{i}}) - \sum_{a=1}^k e_{a}(w_1,\ldots,w_{\mu'_{i+1}})e_{k-a}(w_{\mu'_{i+1}+1}, \ldots, w_{\mu'_i}).
$$
This completes the proof of point (2).
\end{proof}

We finish this section by pointing out that the algebra $\Acr_{GL_\dimv,0}$ was also studied by Cautis and Williams in~\cite{CW18}, who proved the following sharpening of Proposition~\ref{prop-gen}:
\begin{prop}[Lemma 2.18/Corollary 2.21 of~\cite{CW18}]
\label{lem:small-genset}
The classes $\{[\mathcal{Q}_{\varpi_1}^{\otimes m}]\}_{m \in \mathbb{Z}}$, together with the inverse of the class $[\mathcal{O}_{\mathrm{Gr}_{\varpi_{\dimv}}}]$, generate $\Acr_{GL_\dimv,0}$ as an algebra over $\mathbb{Z}[q^{\pm1}]$.
\end{prop}
In particular, it follows that $\Dres(T_\dimv)$ is generated over $\mathbb{Z}[q^{\pm1}]$ by $D_{\varpi^*_\dimv} = (D_1\cdots D_\dimv)^{-1}$ and the first-order $q$-difference operators
\beq
\label{eq:todaclass}
E_{1,(m)}=\sum_{r=1}^\dimv w_r^m \prod_{\substack{s\neq r}} \frac{1}{1-w_{s}w_{r}^{-1}} D_r,\qquad m\in\mathbb{Z}.
\eeq

\subsection{A characterization of the image of $\mathbf{z}^*\hr{\iota_*}^{-1}$ for quiver theories}
Now we return to the setting of general quiver gauge theories  and their Coulomb branch algebras $\Acr_{\Gamma}$. Recall that the homomorphism $ \mathbf{z}^*\hr{\iota_*}^{-1}$ defines an algebra embedding 
$$
\Acr_{\Gamma} \hookrightarrow \Dc_\loc(\Gamma)\simeq \bigotimes_{i\in \Gamma_0^\rf} \Dc_\loc(T_{\dimv_i}) \otimes_{\mathbb{Z}} K_{T_F}(\pt).
$$
The algebra on the right-hand-side contains the subalgebra 
$$
\Acr_{G_V,T_F,0} = \bigotimes_{i\in \Gamma_0^{\rf}} \Dc_{\res}(T_{\dimv_i})\otimes_{\mathbb{Z}} K_{T_F}(\pt).
$$
By Proposition~\ref{prop-gen}, the algebra $\Acr_{\Gamma}$ is generated over $K_{G_V\times T_F\times \widetilde\C^\times}(\pt)$ by the dressed minuscule monoples $E_{\bs n,\bs f}$ and $F_{\bs n,\bs f}$, whose images under $ \mathbf{z}^*\hr{\iota_*}^{-1}$ clearly lie in the subalgebra $\Acr_{G_V,0}\otimes_{\mathbb{Z}}K_{T_F}(\pt)$ by Remark~\ref{rmk:minuscule-res}.
Hence we can regard the localization as defining an algebra embedding
\begin{align}
\label{eq:embed-ch}
\mathbf{z}_{G_V}^* \colon \Acr_{\Gamma} \hookrightarrow \Acr_{G_V,T_F,0},
\end{align}
which can be understood conceptually as coming from the restriction to the zero section of the infinite-rank vector bundle~\eqref{eq:TGN} on $\Gr_{G_V}$.

Now we give a characterization of the image of the algebra embedding~\eqref{eq:embed-ch}. Consider the following non-commutative birational automorphism of the $q$-difference operator algebra $\Dc(\Gamma)$:
\begin{align}
\label{eq:psiGamma}
\mathbf{S}_{N_\Gamma} = \prod_{\chi\in\supp(N_\Gamma)}\Ad_{\Psi(X_\chi^{-1})}.
\end{align}
Note that since all $(T_V\times T_F)$-weight spaces in $N_\Gamma$ have multiplicity one, we can re-express $\mathbf{S}_{N_\Gamma}$ as 
\begin{align}
\label{eq:psiGamma}
\mathbf{S}_{N_\Gamma} = \prod_{a\in\Gamma_1}\Ad_{\mathbf{\Psi}_a},
\end{align}
where the arrow factors are given by
\beq
\label{eq:Psi_a}
\mathbf{\Psi}_{a:i\rightarrow j}= 
\begin{cases}
\prod_{r=1}^{\dimv_i}\prod_{s=1}^{\dimv_j}\Psi\left(u(a) w_{i,r}/w_{j,s}\right) &\text{if} \; a \in \Gamma_1^\rg, \\
\prod_{s=1}^{\dimv_j}\Psi\left(z_i/w_{j,s}\right) &\text{if} \; a \in \Gamma_1^{\rf_+}, \\
\prod_{r=1}^{\dimv_i}\Psi\left(w_{i,r}/z_j\right) &\text{if} \; a \in \Gamma_1^{\rf_-}.
\end{cases}
\eeq

\begin{lemma}
\label{rmk:res-pres}
Suppose $A\in \Dc_\loc(\Gamma)$ is $S_V$-symmetric and satisfies the residue conditions~\eqref{eq:residue-condition} associated to the divisors $d_{\alpha,k}$ for all $\alpha\in\Delta_+(G_V)$ and $k\in\mathbb{Z}$. Then the same is true of $\Ad_{\mathbf{\Psi}_a}(A)$ for any arrow $a\in\Gamma_1$.
\end{lemma}

\begin{proof}
The $S_V$-symmetry of $\Ad_{\mathbf{\Psi}_a}(A)$ is immediate from that of $A$ and $\mathbf{\Psi}_a$. On the other hand, suppose we take 
$$
A = \sum_{\mu} h_\mu D_\mu\in \Dc_\loc(\Gamma),
$$
so that
$$
\Ad_{\mathbf{\Psi}_a}(A) = \sum_{\mu}r_\mu(\bs w,q)h_\mu D_\mu, \qquad r_\mu=\frac{\mathbf{\Psi}_a(\bs w)}{\mathbf{\Psi}_a(q^{2\mu} \bs w)},
$$
where
$$
q^{\la}\bs w = \hr{q^{\ha{\la,\eps_{i,r}}}w_{i,r} ~\Big|~ i \in \Gamma_0^\rg, \, 1 \le r \le \dimv_i}.
$$
Due to the $q$-difference equation~\eqref{eq:qde} we see that $r_\mu(\bs w,q)$ is a rational function. Moreover,
the weight $\alpha$ does not appear in the $(T_V\times T_F)$-subrepresentation $N_a \subset N_\Gamma$ associated to the arrow $a$: this is immediate if $a$ is not a loop based at a gauge node, and in the latter case follows from the fact that for such a loop we must have $u(a)\neq1$. This implies that the function $r_\mu(\bs w,q)$ is regular and nonvanishing at each divisor $d_{\alpha,k}$.
Moreover, for all roots $\alpha = \eps_{i,r} - \eps_{i,s}$ and coweights $\mu$, the $S_V$-symmetry of $\mathbf\Psi_a(\bs w)$ implies that
$$
{r_\mu} \bigg|_{w_{i,s} = q^{2k}w_{i,r}} = {r_{s_\alpha(\mu)+k\alpha}}\bigg|_{w_{i,s} = q^{2k}w_{i,r}},
$$
and the Lemma follows.
\end{proof}

\begin{prop}
\label{prop:conj-crit}
The image of $\Acr_{\Gamma}$ under the algebra embedding~\eqref{eq:embed-ch} admits the following characterization:
\begin{align}
\label{eq:conj-crit}
\Acr_\Gamma = \left\{A\in \Acr_{G_V,T_F,0}~\big|~\mathbf{S}_{N_\Gamma}(A) \in \Acr_{G_V,T_F,0}\right\}.
\end{align}
\end{prop}
\begin{proof}
By Proposition~\ref{prop-gen} we know that $\Acr_\Gamma$ is generated over $K_{G_V\times T_F\times\widetilde{\mathbb{C}}^\times}(\pt)$ by the dressed minuscule monopole operators $E_{\bs n, \bs f},F_{\bs n, \bs f}$. Since $K_{G_V\times T_F\times\widetilde{\mathbb{C}}^\times}(\pt)$ commutes with all $\Psi(X_\chi)$, its elements are invariant under $\mathbf{S}_{N_\Gamma}$ and so clearly satisfy the condition~\eqref{eq:conj-crit}. Moreover, for any coweight $\la$ we see from the $q$-difference equation~\eqref{eq:qde} that conjugating $D_\la$ by $\Psi(X_\chi)$ introduces poles if and only if $\ha{\chi,\la}>0$. So using formula~\eqref{eq:compressed} for the dressed minuscule monopoles, it follows from Lemma~\ref{rmk:res-pres} that $\mathbf{S}_{N_\Gamma}(E_{\bs n, \bs f})$ and $\mathbf{S}_{N_\Gamma}(F_{\bs n, \bs f})$ also remain elements of $\Acr_{G_V,T_F,0}$. Hence any element of $\Acr_\Gamma$ satisfies the condition~\eqref{eq:conj-crit}.

To finish the proof, write $\Acr_{\mathrm{div}}$ for the subalgebra in $\Acr_{G_V,T_F,0}$ given by the right-hand-side of~\eqref{eq:conj-crit}. As in the proof of  Proposition~\ref{prop:ss25}, this algebra has a filtration indexed by the cone of dominant coweights for $G_V$, and we will again prove that $\Acr_\Gamma$ exhausts $\Acr_{\mathrm{div}}$ by induction on the total height of an element.

Given arbitrary element $A \in \Acr_{\mathrm{div}}$, we take a summand $A_\mu\in \Dc^\mu_{\loc}(\Gamma)$ where the coweight $\mu$ is maximal in dominance order among the set of all such summands.  As observed in the proof of Proposition~\ref{prop:ss25}, the summand $A_\mu$ has the form
 $$
A_\mu=  \frac{n_\mu(\bs w,\bs z,\bs u)}{d_\mu(\bs w)}D_\mu
 $$
 for some $n_\mu(\bs w,\bs z,\bs u)\in K_{T_V\times T_F\times\widetilde{\mathbb{C}}^\times}(\pt)$, and the denominator $d_\mu$ is defined by~\eqref{eq:denom}. Since $\mathbf{S}_{N_\Gamma}(A)\in \Acr_{G_V,T_F,0}$, it follows that for each $\chi \in \supp(N_\Gamma)$, conjugating $A_\mu$ by $\Psi(X_\chi^{-1})$ introduces no additional poles to $A_\mu$. This implies that $n_\mu$ is divisible by 
 $$
 \prod_{\substack{\chi \in \supp(N_\Gamma) \\ \ha{\chi,\mu}=-m_\chi<0}}\prod_{k=1}^{m_\chi}(1-q^{2k+1}X_\chi^{-1})
 $$
 in  $K_{T_V\times T_F\times\widetilde{\mathbb{C}}^\times}(\pt)$.
 Hence there exists $f_\mu(\bs w,\bs z,\bs u)\in K_{T_V\times T_F\times\widetilde{\mathbb{C}}^\times}(\pt)$ such that $A_\mu$ has the form
$$
A_\mu = \Ad_{\mathbf{\Psi}^{-1}_{N_\Gamma;\mu}}\circ \Ad_{\mathbf{\Psi}^{-1}_{G_V;\mu}}(f_\mu(\bs w,\bs z,\bs u)D_{\mu}).
$$
Now let $\mathfrak{C}_\mu$ be the refined Weyl chamber~\eqref{eq:ref-Weyl-chamber} in which $\mu$ lies, and $\hc{\omega_1,\ldots,\omega_l}$ the subset of dominant minuscule coweights which generate its cone of integral points. For each such generator $\omega_j$ we have
$$
\prod_{\substack{\chi\in \supp(N_\Gamma) \\ \ha{\chi,\omega_j}\le0}} \Psi(X_\chi^{-1}) =  \prod_{\substack{\chi\in \supp(N_\Gamma) \\ \ha{\chi,\mu}\le0}}\Psi(X_\chi^{-1}),
$$
so recalling formula~\eqref{eq:compressed}, it follows as in the proof of Proposition~\ref{prop:ss25} that there is a sum of products of dressed minuscule monopole operators labelled by the coweights $\hc{\omega_1,\ldots,\omega_l}$ whose degree $\mu$ piece is $A_\mu$. Subtracting this expression from $A_\mu$, we obtain an element of $\Acr_{\mathrm{div}}$ with strictly lower total height, and so by induction it follows that $\Acr_\Gamma = \Acr_{\mathrm{div}}$.
\end{proof}

\begin{remark}
\label{rmk:coprime}
Since we work with the full torus of flavor symmetries, $T_V\times T_F$, the weight multiplicities of $N_\Gamma$ are all 1. So we can re-interpret $\Acr_\Gamma$ as the subalgebra of elements in $\Acr_{G_V,T_F,0}$ which remain there under conjugation by each factor $\mathbf{\Psi}_a$:
$$
\Acr_\Gamma = \bigcap_{a \in \Gamma_1} \Acr_{\Gamma,a},
$$
where
$$
\Acr_{\Gamma,a} = \Acr_{G_V,T_F,N_a} = \hc{A\in \Acr_{G_V,T_F,0}~\big|\,\Ad_{\mathbf{\Psi}_a}(A) \in \Acr_{G_V,T_F,0}}.
$$
\end{remark}

\begin{remark}
D.\,Klyuev pointed out that there should be a more conceptual proof of Proposition~\ref{prop:conj-crit} via a comparison of the characterization of $\Acr_\Gamma$ from the Proposition with a similar one proposed by Teleman in~\cite{Tel21}.
\end{remark}

\section{Quantum cluster algebras}
\label{sec:cluster}

In this section we recall the setup of quantum cluster algebras in a form that is convenient for our purpose in this article. As in~\cite{GS18}, all our cluster seeds will be defined as subsets of the same ambient lattice. Fixing an initial seed $\Pi_0$ with quantum torus $\Tbb_{\Pi_0}$, the quantum torus $\Tbb_{\Pi}$ associated to another seed $\Pi$ is then defined as the image of $\Tbb_{\Pi_0}$ under the corresponding cluster transformation.  The main difference with the standard exposition is that our seeds will only contain the data of the set of \emph{mutable} directions, i.e.\ we do not incorporate a particular choice of frozen basis vectors in the definition of a seed, for the reason that the universally Laurent (also known as the upper cluster $\Xc$-) algebra does not depend on such a choice.

\subsection{Seeds, mutations, and universally Laurent algebras.}

Let us fix from the outset a lattice $\Lambda$ equipped with a skew-symmetric $\mathbb{Z}$-valued bilinear form $(\ast , \ast)$. Associated to this data is an algebra $\Tbb_\Lambda$ 
 over $\mathbb{Z}[q^{\pm 1}]$ called the \emph{quantum torus}, which is spanned over $\mathbb{Z}[q^{\pm 1}]$ by the elements $\hc{Y_\la}_{\la\in \Lambda}$ with the multiplication rule given by
\begin{equation}
\label{quantum.relations}
(-q)^{(\la, \mu)}Y_\la Y_\mu =  Y_{\la+\mu}.
\end{equation}
We caution the reader about the sign in~\eqref{quantum.relations}. In particular, the `classical limit' of our quantum tori will be at $q=-1$ (as is common in skein theory.)

Next, we fix a primitive sublattice $\Lambda_\mut \subset \Lambda$ of rank $r$, where by \emph{primitive} we mean that the abelian group $\Lambda/\Lambda_\mut$ is torsion-free, so that a basis for $\Lambda_\mut$ can be extended to one of $\Lambda$.
For the purposes of this article, we define a \emph{cluster seed} in $\Lambda$ to be a basis
$\Pi$ of the sublattice $\Lambda_\mut$. We refer to the rank $r$ of $\La_\mut$ as the \emph{mutable rank} of the seed. It can often be convenient when writing formulas to label the elements of $\Pi$ by the set $\{1,\ldots, r\}$, although we do not consider such a labelling to be part of our definition of a cluster seed. A cluster seed $\Pi$ defines a collection $\hc{Y_e \,|\, e \in \Pi}$ of (mutable) \emph{cluster $\Xc$-variables}, and a pair of cones in $\La_\mut$:
\[
\Lambda_{\mut}^+({\Pi})=\bigoplus_{e\in \Pi} \mathbb{Z}_{\geq 0} e, \qquad \Lambda_{\mut}^-({\Pi}) = \bigoplus_{e\in \Pi} \mathbb{Z}_{\leq 0} e.
\]
Given a cluster seed $\quiver=\{e_i\}_{i=1}^r$, the skew form on $\Lambda_\mut$ is obviously determined by its values $(e_i,e_j)$ on all pairs of vectors in $\Pi$. It is convenient to encode these values graphically by drawing a directed graph, {known as the \emph{cluster quiver},} with elements of $\Pi$ as vertices, no loops or oriented 2-cycles, and arrows determined by the formula
$$
(e_i,e_j) = \#\hc{i\rightarrow j} -  \#\hc{i\leftarrow j}.
$$

\begin{remark}
\label{rem:Z-valued}
There is a more general setup (required, for example, in order to treat cluster structures associated to moduli spaces of framed $G$-local systems) in which the form $(\ast,\ast)$ is $\Z$-valued on $\La_\mut \times \La$ and $\mathbb{Q}$-valued on $\La \times \La$. In that case the quantum torus is defined over a larger ring $\Z[q^{\pm1/N}]$ for some integer $N$. Since the algebras considered in this article fit into the more restrictive setup, we assume that $(\ast,\ast)$ is $\Z$-valued on the whole of $\La$ in order to simplify the exposition. 
\end{remark}

Suppose that  $\Pi$ is a cluster seed in $\Lambda$, $e_k\in \Pi$, and $\eps\in\{\pm1\}$. Then we define a new cluster seed $\mu^{\eps}_k(\Pi)$ in the same lattice $\Lambda$, called the \emph{mutation of $\Pi$ in direction $e_k$ with sign $\eps$} to be the collection of vectors $\{\mu^{\eps}_k(e_i)\}$ given by
\begin{align}
\label{basis.mutatation}
\mu^{\eps}_k(e_i) = \begin{cases}
 -e_k, &\quad i = k, \\
 e_i + \max\{0,(e_i,e_k)\eps\}e_k, &\quad i\neq k.
\end{cases}
\end{align}
Note that the (unordered) set of vectors $\mu^{\eps}_k(\Pi)$ does not depend on our ordering of the set $\Pi$ used to write the formula~\eqref{basis.mutatation}. Now fix an initial seed $\Pi_0=\{e_\ell\}$, and fix for a concreteness an ordering of $\Pi_0$. Then for any sequence of indices $(k_1,\ldots, k_m)$ with $1\leq k_i\leq r$, the sign coherence of $c$-vectors, cf.~\cite{DWZ07}, states that there exists a unique sequence of signs $(\eps_1,\ldots, \eps_m)$  such that in the cluster seeds \begin{equation}
\label{mutation.quiver}
\Pi_0\stackrel{\mu^{\eps_1}_{k_1}}{\lra} \Pi_1\stackrel{\mu^{\eps_2}_{k_2}}{\lra}\ldots \stackrel{\mu^{\eps_m}_{k_m}}{\lra} \Pi_{m} =\Pi, \qquad  \Pi_j=\hc{e_{\ell}^{(j)}}_{1 \le \ell \le r},
\end{equation}
we have 
$$
f_j=\eps_je_{k_j}^{(j)} \in \Lambda_{\mut}^-(\Pi_0).
$$
In particular, each basis element $\{e_{\ell}^{(j)}\}$ encountered in the sequence~\eqref{mutation.quiver} lies either in $\Lambda_{\mut}^+({\Pi_0})$ or in the opposite cone $\Lambda_{\mut}^-({\Pi_0})$.

If $\quiver$ is another cluster seed in the same lattice $\Lambda$, we say that $\quiver$ is \emph{mutation-equivalent to} $\quiver_0$ if it fits into a chain of sign-coherent mutations~\eqref{mutation.quiver} for some sequence of indices $(k_1, \ldots, k_m)$. We write $\hm{\Pi_0}$ for the set of seeds mutation-equivalent to $\Pi_0$.

Denote by $\Frac(\Tbb_\Lambda)$ the non-commutative field of fractions of $\Tbb_\Lambda$, cf.~\cite[Appendix]{BZ05}.
The {negative} cone $\Lambda^-_\mut(\Pi_0)$ determines a formal completion of the algebra $\Tbb_\Lambda$, which contains the group of non-commutative formal power series with leading term 1:
\[
\widehat{\mathcal{R}}_{\Pi_0}=\Big\{\sum_{\la\in \Lambda_{\mut}^-(\Pi_0)} a_\la(q) Y_\la ~\Big|~ a_0(q)=1, ~~a_\la(q)\in \mathbb{Q}(q)\Big\}.
\]
Now recall the quantum dilogarithm power series $\Psi(X)$ from~\eqref{eq:qdl-def}. If $\Pi$ is a seed mutation equivalent to $\Pi_0$, we use it to define the formal power series 
\begin{align}
\label{eq:qdls}
\Psi_{\Pi;\Pi_0}=  \Psi(Y_{f_m})^{\varepsilon_m}\cdots\Psi(Y_{f_2})^{\varepsilon_2}\Psi(Y_{f_1})^{\varepsilon_1} \; \in \; \widehat{\mathcal{R}}_{\Pi_0}.
\end{align}
We have the following fundamental result, which guarantees that the series $\Psi_{\Pi;\Pi_0}$ is a well-defined function of the set $\Pi$:
\begin{theorem}[{\cite[Theorem 4.1]{Kel17}}]
\label{thm:keller}
 The power series $\Psi_{\Pi;\Pi_0}$ only depends on the basis $\Pi$, and not on the sequence of mutations taking $\Pi_0$ to $\Pi$. 
\end{theorem}

Associated with each seed $\Pi\in |\Pi_0|$ mutation equivalent to $\Pi_0$ is another algebra $\Tbb_{\Pi}$
$$
\Tbb_{\Pi} = {\rm Ad}_{\Psi_{\Pi;\Pi_0}} (\Tbb_\La) \subset \Frac(\Tbb_\Lambda),
$$
Although $\Tbb_{\Pi}$ is clearly isomorphic to the quantum torus algebra $\Tbb_{\Pi_0} = \Tbb_\La$, it is embedded differently into $ {\Frac}(\Tbb_\Lambda)$. The isomorphism of fraction fields
$$
{\rm Ad}_{\Psi_{\Pi;\Pi_0}}\colon{\Frac}(\Tbb_{\Pi_0}) \rightarrow {\Frac}(\Tbb_{\Pi})
$$
is called a \emph{quantum cluster transformation}. Writing it in coordinates relative to the bases $\Pi$ and $\Pi_0$ recovers the original definition of a quantum cluster transformation as given in~\cite{FG09b}. The {\it quantum universally Laurent algebra} is defined to be the intersection
\begin{equation}
\label{quantum.cluster.a}
\mathbb{L}_{\Pi_0}=\bigcap_{\Pi\in |\Pi_0|} \Tbb_{\Pi} \subset {\Frac}(\Tbb_\Lambda).
\end{equation}
The quantum cluster transformation ${\rm Ad}_{\Psi_{\Pi;\Pi_0}}$ defines a canonical isomorphism
$$
\mathbb{L}_{\Pi_0} \simeq \mathbb{L}_{\Pi}.
$$
Another fundamental result, a consequence of ~\cite[Theorem 5.1]{BZ05}, is that the right-hand-side of~\eqref{quantum.cluster.a} can actually be replaced by the finite intersection
\begin{equation}
\label{quantum.cluster.ub}
\mathbb{L}_{\Pi_0}=\Tbb_{\Pi_0} \cap \bigcap_{k=1}^r \Tbb_{\Pi_k} ,\qquad \Pi_k=\mu^{+}_{k}(\Pi_0).
\end{equation}
of quantum tori associated to the seeds reachable from $\Pi_0$ by at most a single mutation.

If $\Pi$ is a cluster seed and $\Pi'\subset \Pi$, we say that $\Pi'$ is obtained from $\Pi$ by \emph{freezing} the directions in $\Pi\setminus\Pi'$. Clearly we have
$$
\Pi'\subset \Pi \implies \Lbb_{\Pi}\subset\Lbb_{\Pi'} 
$$
whenever $\Pi'$ is obtained from $\Pi$ by freezing some subset of vectors. Since $\Lbb_{\hc{e_k}} = \Tbb_{\Pi_0} \cap \Tbb_{\Pi_k}$, we can re-express~\eqref{quantum.cluster.ub} as the statement that
$$
\mathbb{L}_{\Pi}=\bigcap_{k=1}^r \Lbb_{\{e_k\}}, \qquad \Pi = \{e_1,\ldots, e_r\}.
$$
Another useful operation on cluster seeds consists of embedding $\Lambda$ into a larger lattice, while keeping the set $\Pi$ fixed. More precisely, suppose $\widetilde\Lambda$ is a lattice with skew-bilinear form, containing an isometrically embedded copy of  $\Lambda$, and admitting a decomposition $\widetilde\Lambda=\Lambda\oplus\Lambda'$, such that $\Lambda'$ is skew orthogonal to any vector in $\Pi$ (although not necessarily orthogonal to all of $\Lambda$). Then writing $\widetilde\Pi$ for the seed in $\widetilde\Lambda$ given by the same collection of vectors $\Pi$, we have an isomorphism of $\Z[q^{\pm1}]$-modules with the obvious cross-relations between the two factors
\begin{align}
\label{eq:seed-split}
\Lbb_{\widetilde\Pi} \simeq \Lbb_{\Pi} \otimes_{\mathbb{Z}[q^{\pm1}]} \Tc_{\Lambda'}.
\end{align}

We say that two seeds $\Pi',\Pi''\in |\Pi_0|$ mutation-equivalent to $\Pi_0$ are related by a \emph{quasi-permutation} if there is an isometry $\varsigma \colon \Lambda\rightarrow \Lambda$ such that
$$
\varsigma(\Pi') = \Pi''.
$$
The isomorphism of quantum tori induced by $\varsigma$ restricts to an isomorphism of universally Laurent algebras, which we denote by the same symbol:
$$
\varsigma \colon \mathbb{L}^q_{\Pi'} \simeq \mathbb{L}^q_{\Pi''}.
$$
If $\Pi$ is a seed mutation equivalent to $\Pi_0$ and $\varsigma$ is a quasi-permutation such that $\varsigma(\Pi)=\Pi_0$, the corresponding \emph{quantum cluster automorphism} of $\Lbb_{\Pi_0}$ is the algebra automorphism
$$
\varsigma\circ{\rm Ad}_{\Psi_{\Pi;\Pi_0}} \colon \Lbb_{\Pi_0}\rightarrow\Lbb_{\Pi_0}.
$$ The group of all quantum cluster automorphisms under composition is called the \emph{cluster modular group.}

\subsection{Compatible pairs}
We define a \emph{compatible pair} in a lattice $\Lambda$ of rank $d$ to be a pair $(\Pi,\Pi^\vee)$, where $\Pi$ is a cluster seed in $\Lambda$ and $\Pi^\vee=\{\xi_j\}_{j=1}^d$ is an ordered basis for $\Lambda$ for which there exists a (necessarily unique) ordering $\{e_i\}_{i=1}^r$ of $\Pi$ such that
\begin{align}
\label{eq:cpdef}
(e_i,\xi_j)= \delta_{ij}, \quad 1\leq i \leq r, \; 1\leq j\leq d.
\end{align}
{We say that the vectors $\{\xi_i\}$ are \emph{mutable} if $1 \le i \le r$, and are \emph{frozen} otherwise.}
A compatible pair determines a $d\times r$ integer-valued \emph{ensemble matrix} which records how to expand the vectors from the seed $\Pi$ into the basis $\Pi^\vee$:
\beq
\label{eq:ensemble-mat-def}
e_j = \sum_{i=1}^d b_{ij}\eff_i, \qquad 1 \leq j \leq r.
\eeq
Note that if $e_i,e_j\in\Pi$  we have $b_{ij}=(e_i,e_j)$.
We use this matrix to define the rule for signed mutation of compatible pairs in direction $e_k\in\Pi$:
\beq
\label{eq:f-mut}
\mu^\eps_k(\eff_i) =
\begin{cases}
-\eff_k +\sum_{j\neq k} \max\{0,\eps b_{jk}\} \eff_j&\text{if} \; i=k, \\
\eff_i  &\text{if} \; i \ne k,
\end{cases}
\eeq
so that $(\mu_k^\eps(\Pi),\mu_k^\eps(\Pi^\vee))$ defines another compatible pair in $\Lambda$.

The generators $Y_{\xi_j}$ of the quantum torus $\Tbb_\Lambda$ are called the \emph{cluster $\Ac$-variables} associated to the compatible pair $(\Pi,\Pi^\vee)$.
Since $(e_i,\xi_j)\geq0$ for all $i$, it follows from the characterization~\eqref{quantum.cluster.ub} that all cluster $\Ac$-variables lie in the universally Laurent algebra $\Lbb_\Pi$. The \emph{quantum cluster algebra} $\Abb_{\Pi,\Pi^\vee}$ associated to the compatible pair is the algebra generated by the quantum cluster variables from all compatible pairs mutation equivalent to the initial one $(\Pi,\Pi^\vee)$. As just observed, we always have $\Abb_{\Pi,\Pi^\vee}\subset \Lbb_\Pi$, but the algebra $\Abb_{\Pi,\Pi^\vee}$ may in general be smaller.

\begin{remark}
It is sometimes convenient (for example, in the context of moduli spaces of framed local systems) to allow $\Pi^\vee$ to be a basis of the vector space $\La \otimes_\Z \Q$, rather than of the lattice $\La$, satisfying the condition~\eqref{eq:cpdef}.  In this setup one works with a pair of lattices $(\La,\La^\vee)$ with $\La^\vee = \Z\Pi^\vee$. Since this variation is not needed for the class of compatible pairs considered in this article, we stick to the simpler formulation given above.
\end{remark}

\section{Pure gauge theories and the algebraic Whittaker transform}
In this section we recall the cluster description of the quantum Coulomb branch ring $\Acr_{GL_\dimv,0}$ discovered by Cautis and Williams in~\cite{CW18}, and re-interpret aspects of their construction in a way that allows for a generalization to the setting of quiver gauge theories.

\subsection{The Toda chain cluster algebra}
\label{subsec:ltoda}
The underlying quantum cluster algebra for the pure $GL_\dimv$ gauge theory is also known as the $GL_\dimv$ Toda chain cluster algebra. To define it, let $\Prm \simeq \mathbb{Z}^\dimv$ be the weight lattice for $GL_\dimv$, with standard basis $\{\eps_i\}_{i=0}^{d}$ and basis of fundamental weights
$$
\varpi_i = \eps_1 + \ldots + \eps_i, \qquad 1\leq i\leq \dimv.
$$ 
We write $\langle\cdot,\cdot\rangle$ for the standard Euclidean form on $\Prm$ for which $\{\eps_i\}_{i=1}^{\dimv}$ is an orthonormal basis, and
consider the double lattice 
$$
\La = \La_{GL_\dimv} \simeq \Prm \oplus \Prm
$$
equipped with the skew-form
\begin{align}
\label{eq:double-form}
((\lambda_1,\lambda_2),(\mu_1,\mu_2))_{\Lambda} = \langle\lambda_1,\mu_2\rangle - \langle\lambda_2,\mu_1\rangle.
\end{align}
For $\la \in \Prm$ we set
$$
p_\la = (\la,0), \qquad x_\la = (0,\la),
$$
and write $P_\lambda$, $X_\lambda$ for the corresponding generators of the quantum torus $\Tbb_{GL_\dimv,0} = \Tbb_{\La}$ associated to the double lattice. Then the nontrivial commutation relations in $\Tbb_{GL_\dimv,0}$ take the form
$$
P_\lambda X_{\mu} = q^{-2\langle\lambda,\mu\rangle}X_{\mu}P_\lambda, \quad \lambda,\mu\in \Prm.
$$
We will also make the abbreviations
$$
P_i = P_{\eps_i}, \quad X_i = X_{\eps_i} \quad\text{for}\quad i=1,\ldots, d.
$$
Consider the following sets of vectors in $\Lambda$:
\begin{align}
\label{eq:toda-xi-basis}
\Pi &= \{e_{s_i},e_{t_i}\}_{i=1}^{\dimv-1}, & &\hspace{-2cm} e_{s_{i}} =  x_{-\alpha_i}, \quad e_{t_i} = p_{\alpha_i} + x_{\alpha_i},\\
\label{eq:toda-e-basis}
\Pi^\vee &= \{\xi_{s_i},\xi_{t_i}\}_{i=1}^{\dimv}, & &\hspace{-2cm} \xi_{s_{i}} = p_{\varpi_i}+x_{\varpi_i}, \quad \xi_{t_{i}} =  x_{\varpi_i}, 
\end{align}
where $\alpha_i = \eps_i-\eps_{i+1}$ is the $i$-th simple root. Then $(\Pi,\Pi^\vee)$ defines a compatible pair in~$\Lambda$ {with frozen vectors $\xi_{s_\dimv}$ and $\xi_{t_\dimv}$.} We write $\Lbb_{GL_\dimv,0} = \Lbb_{\Pi}$ and $\Abb_{GL_\dimv,0} = \Abb_{\Pi,\Pi^\vee}$ for the universally Laurent algebra defined by the seed $\Pi$ and quantum cluster algebra corresponding to the compatible pair $(\Pi,\Pi^\vee)$ respectively. There is an element $\tau$ of the cluster modular group given by mutating in all directions $\{s_j\}_{j=1}^{\dimv-1}$. 
The quantum cluster automorphism implementing its action on the universally Laurent algebra is 
\begin{align}
\label{eq:tau}
\tau = \varsigma\circ \mathrm{Ad}\left(\prod_{j=1}^{\dimv-1}\Psi(X_{-\alpha_j})^{-1}\right), 
\end{align}
where $\varsigma$ is a generalized permutation whose action on the generators of $\Tbb_{GL_\dimv,0}$ is given by
$$
\varsigma(P_i)=P_i, \quad \varsigma(X_i) = -qP_iX_i,\qquad i=1,\ldots, d.
$$

\begin{lemma}
\label{lem:gens}
The universally Laurent algebra  $\Lbb_{GL_\dimv,0}$ is generated as a $\mathbb{Z}[q^{\pm1}]$-algebra by the elements {$Y_{-\xi_{s_\dimv}}$, $Y_{-\xi_{t_\dimv}}$, and $\tau^m(Y_{\xi_{t_1}})$ with $1-\dimv \le m \le \dimv$.}
\end{lemma}
\begin{proof}
By the results of Goodearl and Yakimov~(\cite{GY21}, Theorem B), it is known that 
$$
\Abb_{GL_\dimv,0}= \Lbb_{GL_\dimv,0}.
$$
Moreover, by the same Theorem both are isomorphic as $\mathbb{Z}[q^{\pm1}]$-algebras to the localization of the ring $A_{q^2}(\mathfrak{n}\hr{(s_0s_1)^{\dimv}}) \otimes_{\mathbb{Z}[q^{\pm2}]} \mathbb{Z}[q^{\pm1}]$ at a pair of generalized minors, where $A_{q^2}(\mathfrak{n}(w))$ is the dual canonical form of the quantum coordinate ring of the unipotent cell associated to the affine Weyl group element $w$ in the loop group $LSL_2$. The Lemma then follows from the fact that under this isomorphism, the two inverted generalized minors correspond to the frozen variables $Y_{\xi_{s_\dimv}}$, $Y_{\xi_{t_\dimv}}$, while the $\tau^m(Y_{\xi_{t_1}})$ with $1-d\leq m\leq d$ correspond to the dual PBW generators of the quantum coordinate ring.
 \end{proof}

\subsection{A faithful representation of $\Lbb_{GL_\dimv,0}$.}
Let {$\Dcr^{\dimv}$ and $\Dcr^\dimv_{\Q}$ be respectively the spaces of all $\Z[q^{\pm1}]$- and $\Q(q)$-valued} functions on the $GL_\dimv$ weight lattice $\Prm$. When the underlying weight lattice is fixed, we omit the upper index and write $\Dcr = \Dcr^\dimv$. Denote by $\Fcr$ and $\Fcr_\Q$ their subspaces of finitely-supported functions. {Throughout this section unless otherwise specified we work over $\Z[q^{\pm1}]$, however although the constructions extend to $\Q(q)$ in an obvious way.}
The following formulas define faithful representations of the quantum torus $\Tbb_{GL_\dimv,0}$ of the initial seed of the Toda chain cluster algebra on $\Fcr$ and $\Dcr$: 
\beq
\label{eq:toda-torus-F-action}
(P_i f)(\mu) = f(\mu-\eps_i), \qquad (X_i f)(\mu) = q^{2\mu_i+\rho_i} f(\mu),
\eeq
where
$$
\rho = (0,-1,\ldots,1-\dimv).
$$
In particular, for a simple root $\alpha_i = \eps_i - \eps_{i+1}$ we have
\beq
\label{eq:Xalpha-action}
(X_{\alpha_i} f)(\mu) =  q^{2(\mu_i-\mu_{i+1})+1} f(\mu).
\eeq

Let 
\begin{align}
\label{eq:Vdef}
\Vcr=\left\{f \in \Fcr \,\big|\, \mu\notin \Prm^+ \implies f(\mu)=0\right\}
\end{align}
be the $\mathbb{Z}[q^{\pm1}]$-submodule of functions vanishing outside the dominant cone $\Prm^+\subset \Prm$. We identify elements of the cone $\Prm^+=\{\lambda_1\geq\lambda_2\geq\cdots\geq \lambda_\dimv\}$ with integer partitions in the usual way. The module $\Vcr$ has a basis consisting of delta-functions $\{\delta_\lambda\}_{\lambda\in\Prm^+}$ supported at a single partition $\lambda\in \Prm^+$:
$$
\delta_\lambda(\mu) =
\begin{cases} 1 &\text{if} \;\; \lambda=\mu, \\
0 & \text{otherwise.}
\end{cases}
$$

Clearly the space $\Vcr$ is not preserved by the $\Tbb_{GL_\dimv,0}$-action on $\Fcr$, since the action involves all translation operators $P_\lambda$, not just those corresponding to dominant $\lambda$. The following result was proved in~\cite{SS25}; {for the reader's convenience we recall the proof.}
\begin{lemma}
\label{lem:pres}
The $\mathbb{Z}[q^{\pm1}]$-submodule $\Vcr\subset\Fcr$ is preserved by the action of the universally Laurent algebra $\mathbb{L}_{GL_\dimv,0}$, and defines a faithful representation of $\mathbb{L}_{GL_\dimv,0}$.
\end{lemma}
\begin{proof}
The condition for a function $f\in \Fcr$ to lie in $\Vcr$ has the form
$$
\mu_{i}-\mu_{i+1} <0 \implies f(\mu)=0.
$$
Given $A\in\Lbb_{GL_\dimv,0}$, we expand it as \begin{align}
\label{eq:a-decomp}
A = \sum_{\nu\in \Prm}a_\nu(\bs X)P_\nu,
\end{align}
so that
{
$$
(A\cdot f)(\mu) = \sum_{\nu\in \Prm}a_\nu\hr{q^{2\mu+\rho}} f(\mu-\nu),
$$
}
where we set
$$
\bs X = \hr{X_1, \ldots, X_{\dimv}} \qquad\text{and}\qquad q^{2\mu+\rho} = \hr{q^{2\mu_1},q^{2\mu_2-1},\ldots,q^{2\mu_{\dimv}+1-\dimv}}.
$$
The preservation of $\Vcr$ by $\Lbb_{GL_\dimv,0}$ will follow if we can show that for all $\mu,\nu\in \Prm$ such that $\mu\notin \Prm^+$ and $\mu-\nu\in \Prm^+$ we have {$a_\nu\hr{q^{2\mu+\rho}}=0$}.
The automorphism part of the cluster mutation in direction $s_i$ consists of conjugation by $\Psi(Y_{-s_i})=\Psi(X_{\alpha_i})$. Hence if $\langle \alpha_i,\nu\rangle=\nu_{i}-\nu_{i+1}\geq0$ we see that $\Ad_{\Psi(X_{\alpha_i})}(P_\nu)$ remains a Laurent polynomial, while 
$$
\nu_{i}<\nu_{i+1}\implies \Psi(X_{\alpha_i}) P_\nu \Psi(X_{\alpha_i})^{-1} = \left(\prod_{k=1}^{\nu_{i+1}-\nu_{i}}\frac{1}{1-q^{2k-1}X_{\alpha_i}}\right) P_\nu.
$$
 Hence $A$ remains Laurent under mutation at $s_i$ if and only if for each $\nu$ with $\nu_i<\nu_{i+1}$, the coefficient $a_\nu$ in~\eqref{eq:a-decomp} has the form
\beq
\label{eq:fnu-factor}
a_\nu(\bs X) = b_\nu(\bs X)\prod_{k=1}^{\nu_{i+1}-\nu_{i}}\left(1-q^{2k-1}X_{\alpha_i}\right) \qquad\text{with}\qquad b_\nu(\bs X)\in \mathbb{Z}[q^{\pm1}][\bs X^{\pm1}].
\eeq
Now let $\mu,\nu\in \Prm$ be such that $\mu\notin \Prm^+$ and $\mu-\nu\in \Prm^+$. Then there exists $1 \le i \le \dimv$ such that $\mu_{i}-\mu_{i+1}<0$, and the condition $\mu-\nu\in \Prm^+$ says that $\mu_{i}-\mu_{i+1}+\nu_{i+1}-\nu_{i}\geq0$. So we have
$$
\nu_{i}-\nu_{i+1}\leq \mu_{i}-\mu_{i+1}<0
$$
and by~\eqref{eq:fnu-factor} we see that $\prod_{k=1}^{\mu_{i+1}-\mu_{i}}\left({1-q^{2k-1}X_{\alpha_i}}\right)$ divides $a_\nu(\bs X)$. Thus  $a_\nu(\bs X)$ acts via a multiplication operator by a Laurent polynomial in variables $q^{2\mu_i+\rho_i}$ divisible by
$$
\prod_{k=1}^{\mu_{i+1}-\mu_{i}}\left({1-q^{2(\mu_{i}-\mu_{i+1}+k)}}\right)
$$
which implies $a_\nu\hr{q^{2\mu+\rho}}=0$ as claimed. 
\end{proof}

\subsection{Algebraic Whittaker transform.}
\label{subsec:whit}

Consider an involution $\vartheta \colon q \mapsto \overline{q}=q^{-1}$ on $\Z[q^{\pm1}]$, and extend it to an 
an involution on $\Dcr$ and $\Fcr$:
\beq
\label{eq:invv}
\vartheta \colon f \mapsto \overline{f}  = \vartheta \circ f.
\eeq
Then we have a {non-degenerate $\Z[q^{\pm1}]$-valued} pairing
\beq
\label{eq:Mpair}
(\cdot, \cdot) \colon \Fcr \otimes_{\Z[q^{\pm1}]} \Dcr_{}\longrightarrow \Z[q^{\pm1}], \qquad (f,\phi)= \sum_{\la\in\Prm}f(\la)\overline{\phi(\la)},
\eeq
which is $\Z[q^{\pm1}]$-linear in the first factor and anti-linear with respect to the involution in the second.
{We can use this pairing to define the action of an operator $A^*$ adjoint to $A \in \Tbb_{GL_\dimv,0}$ on the space $\Dcr$ via
$$
(Af,\phi) = (f,A^*\phi).
$$
In particular, for any $\phi \in \Dcr$ we have
\beq
\label{eq:toda-torus-F-action-adjoint}
(P_i^* \phi)(\mu) = \phi(\mu+\eps_i), \qquad (X_i^* \phi)(\mu) = q^{-2\mu_i-\rho_i} \phi(\mu).
\eeq
}

Now let us consider the ring $\Kc_{\bs w}$ consisting of formal infinite $\Z[q^{\pm1}]$-linear combinations 
\beq
\label{eq:completion}
\psi(\bs w) = \sum_{I\in S}c_I \bs w^I, \quad c_I\in \Z[q^{\pm1}]
\eeq
where
$$
\bs w^I = w_1^{i_1} \ldots w_\dimv^{i_\dimv} \qquad\text{for}\qquad I = (i_1, \ldots, i_\dimv),
$$
and the sum in~\eqref{eq:completion} is taken over a well-ordered set $S \subset \mathbb{Z}^\dimv$ with respect to the reverse-lexicographic order on $\Z^\dimv$ -- that is, there is a minimal element $I \in S$. The action~\eqref{eq:toda-torus-F-action} extends to one on the space
$$
\Dcr_{\bs w} = \Dcr \otimes_{\Z[q^{\pm1}]} \Kc_{\bs w}
$$
of all, not necessarily finitely supported, $\Kc_{\bs w}$-valued functions on $\Z^\dimv$. Declaring $\overline{w_i}= w_i$ for all $1\leq i\leq \dimv$ we extend the pairing~\eqref{eq:Mpair} to 
\beq
\label{eq:Mpair-w}
(\cdot, \cdot) \colon \Fcr \otimes_{\Z[q^{\pm1}]} \Dcr_{\bs w}\longrightarrow \Kc_{\bs w}, \qquad (f,\psi)= \sum_{\la\in\Prm}f(\la)\overline{\psi_\la(\bs w)},
\eeq
where we write $\psi_\la(\bs w)$ for the value of the element $\psi(\bs w) \in \Dcr_{\bs w}$ on the weight $\la \in \Prm$.

In~\cite{SS25}, we construct a distribution $W^{(\dimv)}(\bs w) = W(\bs w) \in \Dcr_{\bs w}$ which we call the \emph{Whittaker distribution} . The value $W_\la(\bs w)$ of the distribution $W(\bs w)$ on the weight $\la \in \Prm$ takes the form of a series
\beq
\label{eq:W-la-series}
W_\la(\bs w) = \bs{w}^\la+\sum_{\mu<\lambda} c_\mu(q)\bs{w}^{\mu},
\eeq
where each monomial has the same total degree in $\bs w$, and the coefficient of the leading term with respect to the dominance order on $\Prm$ is $1$. In particular, it follows from this triangularity property that no nonzero element of $\Tbb_{GL_\dimv,0}$ can annihilate the distribution $W(\bs w)$.

The series $W_{\la}(\bs w)$ are constructed to satisfy \emph{Pieri rules} which can be written in terms of commuting Hamiltonians of the $q$-difference open Toda chain
\beq
H_n = \sum_{|I|=n} P_{\sum_{i\in I}\eps_i} \prod_{m\in G(I)} (1-qX_{-\alpha_m}) \in \Tbb_{GL_\dimv,0},
\eeq
where for a $k$-subset $I\subseteq\{1,\ldots d\}$ we write 
$$
G(I) = \left\{ 1\leq m\leq \dimv~|~ m\notin I~\text{and } m+1 \in I \right\}.
$$
Using the action~\eqref{eq:toda-torus-F-action} of $\Tbb_{GL_\dimv,0}$ on $\Dcr_{\bs w}$ we can write the Pieri rules for $W_\la(\bs w)$ as
\beq
\label{eq:chi-eig}
H_n^* \cdot W_{\la}(\bs w) = e_n(\bs w)W_{\la}(\bs w),
\eeq
where $e_n(\bs w)$ is the $k$-th elementary symmetric function of $\bs w$, i.e.\ the character of the $n$-th exterior power of the defining representation of $GL_\dimv$, and the operator $H_n^*$ is adjoint to $H_n$ with respect to the pairing~\eqref{eq:Mpair-w}.

\begin{example}
For $\dimv = 2$, the Pieri rules read
\begin{align*}
(w_1+w_2)W_{\la_1,\la_2} &= W_{\la_1+1,\la_2} + \hr{1-q^{2(\la_1-\la_2)}} W_{\la_1,\la_2+1}, \\
w_1w_2W_{\la_1,\la_2} &= W_{\la_1+1,\la_2+1}.
\end{align*}
\end{example}

For a dominant weight $\la \in \Prm^+$, the series~\eqref{eq:W-la-series} truncates to a symmetric Laurent polynomial in $\bs w$ known as the \emph{$q$-Whittaker polynomial}
$$
W_{\la}(\bs w) \in K_{GL_{\dimv}\times \widetilde{\mathbb{C}}^\times}(\pt) \simeq \mathbb{Z}[q^{\pm1}][\bs w]^{S_\dimv},
$$
which is the $t=0$ specialization of Macdonald $P$-polynomial $P_\la(q^2,t,\bs w)$, i.e.\
$$
W_\la(\bs w) = P_\la(q^2,0,\bs w).
$$
The $q$-Whittaker polynomials form a basis of the space $K_{GL_{\dimv}\times \widetilde{\C}^\times}(\pt)$ over the ring $\Z[q^{\pm1}]$.

\begin{defn}
The \emph{algebraic Whittaker transform} is the following isomorphism of $\Z[q^{\pm1}]$-modules:
\beq
\label{eq:alg-whit}
\Wc \colon \Vcr \longra  K_{GL_\dimv \times \widetilde\C^\times}(\pt), \qquad f \longmapsto \sum_{\lambda \in \Prm} f(\la)\overline{W_\la}.
\eeq
\end{defn}

Note that the right-hand-side of~\eqref{eq:alg-whit} indeed lies in $K_{GL_\dimv \times \widetilde{\mathbb{C}}^\times}(\pt)$ since elements of $\Vcr$ are non-vanishing only on dominant $\la$, in which case $W_\la(\bs w)$ is a $q$-Whittaker polynomial. That the resulting map $\Wc$ is indeed an isomorphism follows since it takes the delta-function basis $\{\delta_\la\}_{\la\in P^+}$ in $\Vcr$ to the basis $\{\overline{W(\la)}\}_{\la\in\Prm^+}$ in $K_{GL_\dimv\times \widetilde{\mathbb{C}}^\times}(\pt)$.

\subsection{Isomorphism with the quantized Coulomb branch ring}

Being specializations of Macdonald polynomials, the $q$-Whittaker polynomials satisfy a dual set of difference equations with respect to the variables $\bs w$: we have
\begin{align}
\label{eq:qwhit-dual-eigen}
(\iota_*)^{-1}([\mathcal{O}_{\Gr_{\varpi_n}}]) \cdot \overline{W}_\la = q^{2\sum_{r=1}^n\la_r}\overline{W_\la},\qquad 1\leq k\leq \dimv,
\end{align}
where 
$$
(\iota_*)^{-1}([\mathcal{O}_{\Gr_{\varpi_n}}]) = \sum_{|J|=n}\prod_{\substack{r \in J \\ s \notin J}}\frac{1}{1-w_sw_r^{-1}} \prod_{j\in J}D_j
$$
is the $q$-difference operator given by the image of the class of the structure sheaf of $\Gr_{\varpi_n}$ under the localization embedding, cf.~\eqref{Eims}. They also satisfy the `raising operator' identities describing the action of the class of the determinant of the tautological quotient bundle
$$
(\iota_*)^{-1}([\det\Qc_{\varpi_n}]) = \sum_{|J|=n}\prod_{\substack{r \in J \\ s \notin J}}\frac{1}{1-w_sw_r^{-1}} \prod_{j\in J}w_jD_j, 
$$
which read
$$
[\det\Qc_{\varpi_n}] \cdot \overline{W}_\la= q^{2\sum_{r=1}^n\la_r}\overline{W}_{\la+\varpi_n},\quad 1\leq n \leq \dimv.
$$
Hence we have the following intertwining relations identifying the action of the frozen $A$-variables from the quantum cluster algebra $\Abb_{GL_\dimv,0}$ with those of rescaled minuscule monopole operators:
\begin{align}
\label{eq:A-inter}
\Wc\circ Y_{\xi_{t_n}} &= q^{-n(n-1)/2}[\Oc_{\Gr_{\varpi_n}}]\circ\Wc,\\
\Wc\circ Y_{\xi_{s_n}} &= {(-q)^n q^{-n(n-1)/2}}[\det\Qc_{\varpi_n}]\circ\Wc.
\end{align}

Before we proceed to the main result of this section, note that it is easy to check directly that the elements $H_n \in \Tbb_{GL_\dimv,0}$ remain Laurent under mutations in directions $s_k$ and $t_k$ for all $1 \le k \le n$, and thus lie in $\Lbb_{GL_\dimv,0}$. We also recall the following result, cf. ~\cite[Lemma 6.23]{SS25} and~\cite[Theorem 3.8]{DFK16}:
\begin{lemma}
\label{lem:time-translation}
For all $m\in\mathbb{Z}$, we have
$$
\hs{H_{\pm1},\tau^m(Y_{\xi_{t_1}})} = \pm(q-q^{-1}) \tau^{m\pm1}(Y_{\xi_{t_1}}), \\
$$
where
$$
H_{-n} = H_\dimv^{-1} H_{d-n}.
$$
\end{lemma}

The next result is obtained in~\cite{CW18}. For the reader's convenience we recall a sketch of its proof in the language introduced above.
\begin{prop}
\label{prop:alg-whit}
There is a unique algebra isomorphism 
\begin{align}
\label{eq:alg-whit-def}
\mathbb{W}\colon \mathbb{L}_{GL_\dimv,0}\longrightarrow \mathscr{A}_{GL_\dimv,0},
\end{align}
such that
\begin{align}
\label{eq:el-inter}
\mathbb{W}(H_n) &= e_n(\bs w), \\
\mathbb{W}(Y_{\xi_{t_n}}) &= q^{-n(n-1)/2}[\Oc_{\Gr_{\varpi_n}}],\\
\mathbb{W}(Y_{\xi_{s_n}}) &= (-q)^n q^{-n(n-1)/2}[\det\Qc_{\varpi_n}].
\end{align}
It satisfies the intertwining relation
\beq
\label{eq:aut-inter}
\mathbb{W}\circ \tau = \gamma \circ \mathbb{W},
\eeq
where $\gamma $ is the automorphism of $\Acr_{GL_\dimv,0}$ defined by
$$
\gamma(w_i)=w_i,\quad \gamma(D_i) = {-q} w_iD_i, \quad i=1,\ldots,d.
$$
\end{prop}

\begin{proof}
Using the isomorphism $\Wc$ we can identify both $\mathscr{A}_{GL_\dimv,0}$ and $\Lbb_{GL_\dimv,0}$ as subrings of the endomorphism ring $\End_{\Z[q^{\pm1}]}(K_{GL_\dimv\times\widetilde\C^\times}(\pt))$, and we claim that these subrings are identical. Indeed, comparing the relations from Lemma~\ref{lem:time-translation} and the evident commutation relations between $e_1(\bs w^{\pm1})$ and $[\Oc_{\Gr_{\varpi_1}}]$, it follows that the action of $\tau^m(Y_{\xi_{t_1}})$ coincides up to multiplication by a power of $q$ with that of the class $[\Qc_{\varpi_1}^{\otimes m}]$. Similarly, the frozen $\Ac$-variables $Y_{\xi_{t_\dimv}} = X_{\varpi_\dimv}$ and $Y_{\xi_{s_\dimv}} = q^\dimv P_{\varpi_\dimv} X_{\varpi_\dimv}$ are intertwined respectively with powers of $q$ times $[\Oc_{\Rc_{\varpi_\dimv}}]$ and $[\det\Qc_{\varpi_\dimv}]$. Since the $\tau^m(Y_{\xi_{t_1}})$, $Y_{-\xi_{t_\dimv}}$,$Y_{-\xi_{s_\dimv}}$ generate $\Lbb_{GL_\dimv,0}$ by Lemma~\ref{lem:gens}, while the $[\Qc_{\varpi_1}^{\otimes m}]$ together with the inverse of $[\Oc_{\Gr_{\varpi_\dimv}}]$ generate $\Acr_{GL_\dimv,0}$ by Lemma~\ref{lem:small-genset} {and $[\det\Qc_{\varpi_\dimv}]^{-1} \in \Acr_{GL_\dimv,0}$}, the claim follows.
\end{proof}

\section{Cluster description of $\Acr_\Gamma$ for $\Gamma$ without 1-cycles}

\subsection{From gauge theory quivers to cluster seeds}
\label{cluster-quiver-sect}
Consider the quantized $K$-theoretic Coulomb branch algebra $\Acr_{\Gamma}$ of a  quiver gauge theory corresponding to a weighted quiver $\Gamma$ without 1-cycles as in the setup of Section~\ref{sec:bfn}.  To such a datum we associate a lattice
\begin{align}
\label{eq:lambda-def}
\Lambda_{|\Gamma|} = \bigoplus_{i\in \Gamma^{\rg}_0}\Lambda_{GL_{\dimv_i}}\oplus\bigoplus_{k\in \Gamma^{\rf}_0}\mathbb{Z}\zeta_k\oplus \bigoplus_{r=1}^g\mathbb{Z}\upsilon_r.
\end{align}
We equip $\Lambda_{|\Gamma|} $ with the skew form such that the summands in~\eqref{eq:lambda-def} are pairwise orthogonal, and its restriction to the $\Lambda_{GL_{\dimv_i}}$ coincides with~\eqref{eq:double-form} while its restriction to all other summands is zero. In other words, the lattice $\Lambda_{|\Gamma|}$ has a Darboux basis 
$$
\left\{p_{i,m},x_{i,m} ~\big|~ i\in \Gamma^{\rg}_0, ~1\leq m\leq \dimv_i\right\} \sqcup \big\{\zeta_{k}~\big|~ k\in \Gamma^{\rf}_0\big\} \sqcup \{ \upsilon_r \}_{r=1}^{g}
$$
where 
$$
(p_{i,m},x_{j,n}) = \delta_{i,j}\delta_{m,n}
$$
and $\zeta_{k}$, $\upsilon_r$ lie in the kernel of the pairing.

As indicated by the notation, the lattice $\Lambda_{|\Gamma|}$ and its skew form do not depend on the orientation of the edges in the weighted quiver $\Gamma$, only on the dimension vector $\bs\dimv = (\dimv_i \,|\, i \in \Gamma_0^\rg)$ and the underlying undirected graph which we denote by $|\Gamma|$. The quantum torus $\Tbb_{\Lambda_{|\Gamma|}}$ associated to the lattice~\eqref{eq:lambda-def} can be written as
\beq
\label{eq:qtor-ident}
\Tbb_{\Lambda_{|\Gamma|}} \simeq \bigotimes_{i \in \Gamma_0^\rg} \Tbb_{GL_{\dimv_i},0} \otimes_{\mathbb{Z}} K_{T_F}(\pt),
\eeq
where the first tensor product is taken over the ring $\Z[q^{\pm1}]$.

\begin{defn}
\label{eq:seed-def}
We define the cluster seed $\quiver_\Gamma$
in the lattice $\Lambda_{|\Gamma|}$ to be the following set of vectors:
\begin{itemize}
\item for each gauge node $i\in \Gamma_0^\rg$ and $1\leq m<\dimv_i$, the vectors
\beq
\label{eq:est}
e_{s_{i,m}} = x_{i,m+1} - x_{i,m},\qquad e_{t_{i,m}} = p_{i,m}-p_{i,m+1} + x_{i,m}-x_{i,m+1};
\eeq
\item for each arrow $a \colon i\rightarrow j\in \Gamma_1$, the vector
\beq
\label{eq:earr}
e_a =
\begin{cases}
p_{j,\dimv_j} - p_{i,1} + \upsilon(a) &\text{if} \quad a \in \Gamma_1^\rg, \\
\zeta_j - p_{i,1} &\text{if} \quad a \in \Gamma_1^{\rf_+}, \\
p_{j,\dimv_j} - \zeta_i &\text{if} \quad a \in \Gamma_1^{\rf_-},
\end{cases}
\eeq
where
$$
\upsilon(a)=\sum_{r=1}^{g}\theta_r(a)\upsilon_r,
$$
and the 1-cochains $\theta_r$ are as defined in Section~\ref{subsec:coulomb-def}.
\end{itemize}
\end{defn}

\begin{remark}
\label{rem:cluster-quiver}
The vectors in $\quiver_\Gamma$ are labelled by the simple roots of the gauge group $G_V$ and the lowest weights of the $G_V \times T_F$-representation $N_\Gamma$. Indeed, to the $m$-th simple root $\alpha_{i,m}$ of the factor $GL(V_i) \subset GL_V$ we associate a pair of vectors $e_{s_{i,m}} = x_{-\alpha_{i,m}}$ and $e_{t_{i,m}} = p_{\alpha_{i,m}}+x_{\alpha_{i,m}}$. To the direct summand $N_a \subset N_\Gamma$ contributed by the arrow $a \in \Gamma_1$, we associate the vector $e_a = p_{\omega}+\upsilon_{\theta}$, where $\omega$ is the lowest weight of $N_a$ with respect to the group $G_V \times T_W\times T_H$, where we understand $p_k = \zeta_k$ when $k \in \Gamma_0^\rf$.
\end{remark}

In order to check that the set $\Pi_\Gamma$ indeed defines a seed we need to verify that it can be extended to a basis for $\Lambda_{|\Gamma|}$. We can give a constructive proof by equipping $\Gamma$ with a `rooting' consisting of the choice of a  \emph{marked} gauge node $i_{C}\in \Gamma_0^\rg$ within each connected component $C\in \pi_0(\Gamma)$ of $\Gamma$. We write $\Gamma^\dag$ for a quiver equipped with such a choice of marked gauge nodes. Associated to such a $\Gamma^\dag$ is an additional collection $\Pi^{\mathrm{fr}}_{\Gamma^\dag}$ of vectors in $\Lambda$ defined as follows:

\begin{itemize}
\item for each gauge node $i\in \Gamma^\rg_0$, we define a frozen basis vector
$$
e_i = x_{i,1};
$$
\item for each connected component $C\in \pi_0(\Gamma)$ with its root gauge node $i_C$, we define a frozen basis vector
\beq
\label{eq:eroot}
e_C = p_{i_C, \dimv_{i_C}}.
\eeq
\end{itemize}

\begin{lemma}
\label{lem:cards}
We have $\mathrm{rank}(\Lambda_{|\Gamma|}) = |\Pi_\Gamma| + |\Gamma_0^\rg| + |\pi_0(\Gamma)|$.
\end{lemma}
\begin{proof}
Writing $\dim(V)$ for the sum of dimensions of the vector spaces $V_i$ at gauge nodes, $\dim(W) = |\Gamma^{\rf}_0|$ for the total framing dimension, and $b_1(\Gamma) = \dim(H^1(\Gamma,\C^\times))$ for the first Betti number of $\Gamma$ we have
\beq
\label{eq:lambda-dim}
\mathrm{rank}(\Lambda_{|\Gamma|}) = 2\dim(V) + \dim(W) + b_1(\Gamma).
\eeq
Since each node of $\Gamma$ is either gauge or flavor and the number of flavor nodes is $\dim(W)$, we see that 
$$
|\Pi_\Gamma| + |\Gamma_0^\rg| + |\pi_0(\Gamma)| = 2\dim(V) + |\Gamma_1| + \dim(W) - |\Gamma_0| + b_0(\Gamma).
$$
Now the result follows from comparing the following two expressions for the Euler characteristic of $\Gamma$:
$$
\chi(\Gamma) = |\Gamma_0| - |\Gamma_1| = b_0(\Gamma) - b_1(\Gamma).
$$
\end{proof}

\begin{lemma}
\label{lem:linalg}
The set of vectors $\quiver_{\Gamma^\dag} = \quiver_\Gamma \cup \Pi^{\mathrm{fr}}_{\Gamma^\dag}$ forms a basis for $\Lambda_\Gamma$.
\end{lemma}

\begin{proof}
By Lemma~\ref{lem:cards} it is enough to prove that the vectors in $\quiver_{\Gamma^\dag}$ span $\Lambda_\Gamma$. Defining $\Lambda'$ as the span of $\quiver_{\Gamma^\dag}$, we see that $\La'$ contains vectors
$$
x_{i,n} = e_i + \sum_{r=1}^{n-1} e_{s_{i,r}}, \qquad p_{i,m}-p_{i,m+1} = e_{s_{i,m}}+e_{t_{i,m}}
$$
for all $i \in \Gamma_0^\rg$, $1 \le n \le \dimv_i$, and $1 \le m < \dimv_i$. So it suffices to show that the lattice $\Lambda'$ contains all vectors $\upsilon_r$, $\zeta_k$, together with the vectors $p_{i,1}$ for each gauge node $i$. Observe that for each $i \in \Gamma_0^\rg$, the vector
$$
\lambda(i) = p_{i,1}-p_{i,\dimv_{i}}= \sum_{m=1}^{\dimv_{i}-1}(e_{s_{i,m}}+e_{t_{i,m}})  
$$
lies in $\Lambda'$. So we can define a $\Lambda'$-valued 1-cochain on $\Gamma$ whose value on each arrow $a \colon i\rightarrow j$ in $\Gamma_1^\rg$ is given by
$$
{\eta(a) = e_a+\lambda(j) = p_{j,\dimv_j}-p_{i,\dimv_i}+\upsilon(a).}
$$
Hence if $(\gamma_1,\ldots, \gamma_g)$ is the homology basis dual to $([\theta_1],\ldots, [\theta_g])$, we have $\eta(\gamma_i) = \upsilon_i$, which shows that each $\upsilon_i$ lies in $\Gamma'$. Similarly, given a general gauge node $i$ lying in connected component $C$ of $\Gamma$ and a path $\gamma$ from the marked gauge node $i_C$ to $i$, we have {$e_C+\eta(\gamma)= p_{i,\dimv_i} + \upsilon(\gamma)$}, which shows $p_{i,\dimv_i} \in \Lambda'$, and hence all $p_{i,m}$ lie in $\Lambda'$.
Finally, using similar paths ending at the framing nodes we get that $\zeta_k\in \Lambda'$ for all $k$, which completes the proof.
\end{proof}

\subsection{Laurent algebra with frozen matter directions}
\label{subsec:Laurent-frozen}

\begin{defn}
We write {$\quiver_{G_V,T_F,0}$} for the seed in the lattice $\Lambda_{|\Gamma|}$ obtained from $\quiver_\Gamma$ by freezing, or simply omitting in the notations of Section~\ref{sec:cluster}, the vectors $e_a$ for all arrows $a \in \Gamma_1$. We denote its universally Laurent algebra by {$\Lbb_{G_V,T_F,0} = \Lbb_{\quiver_{G_V,T_F,0}}$} and refer to it as the \emph{universally Laurent algebra with frozen matter directions} associated to the quiver gauge theory.
\end{defn}

The universally Laurent algebra with frozen matter directions is easy to describe. Indeed, as observed in the proof of Lemma~\ref{lem:linalg}, the following set of vectors forms a basis for $\Lambda_{|\Gamma|}$
$$
\{e_{s_{i,m}},e_{t_{i,m}}~|~ i\in \Gamma_0^\rg,~1\leq m < \dimv_i\} \sqcup \{\zeta_k\}_{k \in \Gamma_0^\rf} \sqcup \{\upsilon_r\}_{r=1}^g.
$$
Hence writing
$$
\Lbb_{G_V,0} = \bigotimes_{i\in \Gamma_0^\rg}\Lbb_{\mathrm{GL}_{\dimv_i},0},
$$ for the tensor product over $\mathbb{Z}[q^{\pm1}]$ of the Toda universally Laurent algebras associated to each gauge node of $\Gamma$,  
by~\eqref{eq:seed-split} we have an isomorphism
\beq
\label{eq:local-laurent-ring}
\Lbb_{G_V,T_F,0}
\simeq \Lbb_{G_V,0} \otimes_\Z K_{T_F}(\pt),
\eeq
with the generators $Y^{\pm1}_{\zeta_k}$, $Y^{\pm1}_{\upsilon_r}$ of $K_{T_F}(\pt)$ being central in $\Lbb_{G_V,T_F,0}$. In particular, $\Lbb_{G_V,T_F,0}$ is generated by its central subalgebra $K_{T_F}(\pt)$ together with the generators from Lemma~\ref{lem:gens} associated to each factor in $G_V$.

Let $\Dcr^\Gamma$ be the space of all $K_{T_F}(\pt)$-valued functions on the weight lattice $\Prm$ of $G_V$, $\Fcr^\Gamma \subset \Dcr^\Gamma$ be the subset of finitely-supported functions, and $\Vcr^\Gamma \subset \Fcr^\Gamma$ be the subspace of finitely-supported functions vanishing outside the dominant cone $\Prm^+$. As before we write $\Dcr^\Gamma_\Q = \Dcr^\Gamma \otimes_{Z[q^{\pm1}]} \Q(q)$ and $\Fcr^\Gamma_\Q = \Fcr^\Gamma \otimes_{Z[q^{\pm1}]} \Q(q)$. Consider the ring 
$$
\Kc^\Gamma_{\bs w} = K_{T_F}(\pt)\bigotimes_{i\in \Gamma^\rg_0} \mathcal{K}_{\bs w^{(i)}}
$$
where $\mathcal{K}_{\bs w^{(i)}}$ is the copy of the ring defined in section~\eqref{subsec:whit}, associated to the $GL_{\dimv_i}$ factor of the gauge group. Then we can define the space
$$
\Dcr^\Gamma_{\bs w} = \Dcr^\Gamma \otimes_{\Z[q^{\pm1}]} \Kc^\Gamma_{\bs w}
$$
and observe that the formula~\eqref{eq:Mpair-w} defines a pairing
\beq
\label{eq:wwpair}
(\cdot, \cdot) \colon \Fcr^\Gamma \otimes_{\Z[q^{\pm1}]} \Dcr^{\Gamma}\longrightarrow \Kc^\Gamma_{\bs w},
\eeq
where we declare all equivariant parameters of $T_F$ to be invariant under the involution~\eqref{eq:invv}. In this space we have the $G_V$-Whittaker distribution 
$$
W^\Gamma(\bs w)  = \prod_{i} W^{\dimv_i}\big(\bs w^{(i)}\big) \in \Dcr^{\Gamma}_{\bs w},
$$
which takes the value
$$
W^\Gamma_{\la}(\bs w) = \prod_{i} W^{\dimv_i}_{\la^{(i)}}\big(\bs w^{(i)}\big)
$$
on the weight $\la = (\la^{(i)} \,|\, i \in \Gamma_0^\rg) \in \Prm$ of the group $G_V$.

We equip $\Fcr^\Gamma$ with an action of the quantum torus $\Tbb_{\Lambda_{|\Gamma|}}$ using the standard formulas~\eqref{eq:toda-torus-F-action} to define the action of the generators $P_{i,m}$, $X_{i,m}$ corresponding to gauge nodes, and letting the generators from $K_{T_F}(\pt)$ act by scalar multiplication.  By Lemma~\ref{lem:pres} the universally Laurent algebra with frozen matter directions $\Lbb_{G_V,T_F,0}$ preserves the subspace $\Vcr^\Gamma\subset\Fcr^\Gamma$ of functions vanishing outside the dominant cone in $G_V$, and as in Proposition~\ref{prop:alg-whit} the $G_V$-Whittaker transform
\beq
\label{eq:alg-whit-gamma}
\Wc^\Gamma \colon \Vcr^\Gamma \simeq  K_{G_V\times T_F\times \widetilde\C^\times}(\pt), \qquad f \longmapsto \sum_{\lambda \in \Prm} f(\la)\overline{W^\Gamma_{\la}}
\eeq
induces an isomorphism of $\Z[q^{\pm1}]$-algebras 
\beq
\label{eq:boldW-Gamma}
\mathbb{W}_{\Gamma}\colon\Lbb_{G_V,T_F,0} \simeq \Acr_{G_V,T_F,0}
\eeq
identifying both sides with subrings in the endomorphism ring of $ K_{G_V\times T_F\times \widetilde\C^\times}(\pt)$.

\subsection{Reversing the orientation of an arrow in $\Gamma$}
\label{subsec:bibax}

Recall that the definition of the lattice $\Lambda_{|\Gamma|}$ and its skew-form in~\eqref{eq:lambda-def} does not depend on the orientation of the weighted quiver $\Gamma$. So if $\Gamma_{\overline{a}}$ is the quiver obtained from $\Gamma$ by replacing a single arrow $a \colon i\rightarrow j \in\Gamma_1$ connecting distinct vertices of $\Gamma$ with an arrow $\overline{a} \colon j\rightarrow i $, we obtain two distinct cluster seeds 
$$
\quiver = \quiver_{\Gamma},\qquad \quiver'= \quiver_{\Gamma_{\overline{a}}}
$$
in the same lattice $\Lambda_{|\Gamma|}$. The corresponding quantum universally Laurent algebras are two subalgebras in the same quantum torus:
\begin{align}
\label{eq:incs}
\Lbb_{\quiver} \hookrightarrow \Tbb_{\Lambda_{|\Gamma|}} \hookleftarrow \Lbb_{\quiver'}.
\end{align}
The cluster seeds $\quiver$ and $\quiver'$ are almost identical. Indeed, to obtain $\quiver'$ from $\quiver$ it suffices to replace the single vector $e_a \in \quiver$ by the vector $e_{\overline a} \in \quiver'$. The two vectors in question take the form
$$
e_ a = p_{j,\dimv_j} - p_{i,1} + \upsilon(a), \qquad e_{\overline{a}} = p_{i,\dimv_i} - p_{j,1} - \upsilon(a)
$$
since $\overline{a} = - a$ as 1-chains. In particular, note that $\quiver_{G_V,0}$ is contained on both $\quiver$ and $\quiver'$, and so the inclusions~\eqref{eq:incs} factor through
\beq
\label{eq:small-incs}
\Lbb_{\quiver} \hookrightarrow \Lbb_{G_V,T_F,0} \hookleftarrow \Lbb_{\quiver'}.
\eeq

In this section we describe a {\emph{Baxter sequence}} of seed mutations,
$$
\bs \mu_{ a} =  \mu_{i_M}\circ\cdots \circ  \mu_{i_1}, \qquad M = (2\dimv_i-1)(2\dimv_j-1)
$$ 
such that 
$$
\quiver_{\Gamma_{\overline a}}= \bs \mu_{a}(\quiver_\Gamma).
$$
These mutations take place entirely within the sub-seed
\begin{align}
\label{eq:arrow-subseeds}
\quiver_a = \mathrm{span}\left\{e_{s_{i,m}},e_{t_{i,m}} \right\}_{m=1}^{\dimv_i-1}\oplus \mathbb{Z}\langle e_ a\rangle \oplus \mathrm{span}\left\{e_{s_{j,n}},e_{t_{j,n}} \right\}_{n=1}^{\dimv_j-1}.
\end{align}
which we represent graphically for $\dimv_i=4$ and $\dimv_j=3$ in Figure~\ref{fig-Qmn}. To describe the mutation sequence $\bs \mu_{a}$, it is convenient to label the vectors in the basis~\eqref{eq:arrow-subseeds} by integers $\{1,\ldots, 2(\dimv_i+\dimv_j)-3\}$ as indicated in Figure~\ref{fig-Qmn}. Precisely, the labeling reads as follows:
\begin{align*}
e_{s_{j,k}} &\mapsto 2k-1, & &e_{t_{j,k}} \mapsto 2k, & &e_a \mapsto 2\dimv_j-1, \\
e_{s_{i,k}} &\mapsto 2(\dimv_j+k)-1, & &e_{t_{i,k}} \mapsto 2(\dimv_j+k-1).
\end{align*}

\begin{figure}[h]
\begin{minipage}{.5\textwidth}
\centering
\begin{tikzpicture}[every node/.style={inner sep=0, minimum size=0.45cm, draw, circle}, x=0.5cm, y=0.5cm]

\def\m{3};
\def\n{2};

\foreach \i in {-\n,...,\m}{
  \pgfmathparse{int(2*\i)}\edef\y{\pgfmathresult};
  \pgfmathparse{int(2*(\n+\i)+1)}\edef\v{\pgfmathresult};
  \node (0_\y) at (0,\y) {\tiny{$\v$}};
}

\pgfmathparse{int(1-\n)}\edef\t{\pgfmathresult};
\foreach \i in {\t,...,\m}{
  \pgfmathparse{int(2*\i-1)}\edef\y{\pgfmathresult};
  \pgfmathparse{int(2*(\n+\i))}\edef\v{\pgfmathresult};
  \node (2_\y) at (2,\y) {\tiny{$\v$}};
}

\pgfmathtruncatemacro{\nminusone}{\n-1}

\foreach \k in {0,...,\nminusone}{
  \pgfmathparse{int(2*(\k-\n))}\edef\y{\pgfmathresult}
  \pgfmathparse{int(\k+1)}\edef\kp{\pgfmathresult}
  \node[anchor=east, xshift=-8pt, yshift=-1pt, draw=none] at (0_\y) {\footnotesize$s_{j,\kp}$};
}

\node[anchor=east, xshift=-8pt, yshift=-1pt, draw=none] at (0_0) {\footnotesize$a$};

\foreach \k in {1,...,\m}{
  \pgfmathparse{int(2*\k)}\edef\y{\pgfmathresult}
  \node[anchor=east, xshift=-8pt, yshift=-1pt, draw=none] at (0_\y) {\footnotesize$s_{i,\k}$};
}

\foreach \k in {1,...,\n}{
  \pgfmathparse{int(2*(\k-\n)-1)}\edef\y{\pgfmathresult}
    \pgfmathparse{int(\k+1)}\edef\kp{\pgfmathresult}
  \node[anchor=west, xshift=8pt, yshift=-1pt, draw=none] at (2_\y) {\footnotesize$t_{j,\k}$};
}

\foreach \k in {1,...,\m}{
  \pgfmathparse{int(2*\k-1)}\edef\y{\pgfmathresult}
  \node[anchor=west, xshift=8pt, yshift=-1pt, draw=none] at (2_\y) {\footnotesize$t_{i,\k}$};
}

\foreach \i in {1,...,\m}{
  \pgfmathparse{int(2*\i)}\edef\y{\pgfmathresult};
  \pgfmathparse{int(2*\i-1)}\edef\t{\pgfmathresult};
  \draw [->] (0_\y.-10) to (2_\t.140);
  \draw [->] (0_\y.-35) to (2_\t.165);
}

\foreach \i in {2,...,\m}{
  \pgfmathparse{int(2*\i)}\edef\y{\pgfmathresult};
  \pgfmathparse{int(2*\i-3)}\edef\t{\pgfmathresult};
  \draw [<-] (0_\y) to (2_\t);
}

\foreach \i in {1,...,\m}{
  \pgfmathparse{int(2*\i-1)}\edef\t{\pgfmathresult};
  \pgfmathparse{int(2*\i-2)}\edef\y{\pgfmathresult};
  \draw [<-] (0_\y) to (2_\t);
}

\draw [->] (0_0) to (0_2);
\draw [->] (0_0) to (0_-2);

\foreach \i in {1,...,\n}{
  \pgfmathparse{int(-2*\i)}\edef\y{\pgfmathresult};
  \pgfmathparse{int(1-2*\i)}\edef\t{\pgfmathresult};
  \draw [->] (0_\y.10) to (2_\t.-140);
  \draw [->] (0_\y.35) to (2_\t.-165);
}

\foreach \i in {2,...,\n}{
  \pgfmathparse{int(-2*\i)}\edef\y{\pgfmathresult};
  \pgfmathparse{int(3-2*\i)}\edef\t{\pgfmathresult};
  \draw [<-] (0_\y) to (2_\t);
}

\foreach \i in {1,...,\n}{
  \pgfmathparse{int(1-2*\i)}\edef\t{\pgfmathresult};
  \pgfmathparse{int(2-2*\i)}\edef\y{\pgfmathresult};
  \draw [<-] (0_\y) to (2_\t);
}

\end{tikzpicture}

\caption{Seed $\quiver_{a:i\rightarrow j} \subset \quiver_\Gamma$ for $\dimv_i=4, \dimv_j=3$.}
\label{fig-Qmn}
\end{minipage}%
\begin{minipage}{.5\textwidth}
\centering

\begin{minipage}{.5\textwidth}
\centering
\begin{tikzpicture}[every node/.style={inner sep=0, minimum size=0.45cm, draw, circle}, x=0.5cm, y=0.5cm]

\def\m{3};
\def\n{2};
\foreach \i in {-\m,...,\n}{
  \pgfmathparse{int(2*\i)}\edef\y{\pgfmathresult};
  \pgfmathparse{int(2*(\m+\i)+1)}\edef\v{\pgfmathresult};
  \node (0_\y) at (0,\y) {\tiny{$\v$}};
}

\pgfmathparse{int(1-\m)}\edef\t{\pgfmathresult};
\foreach \i in {\t,...,\n}{
  \pgfmathparse{int(2*\i-1)}\edef\y{\pgfmathresult};
  \pgfmathparse{int(2*(\m+\i))}\edef\v{\pgfmathresult};
  \node (2_\y) at (2,\y) {\tiny{$\v$}};
}

\pgfmathtruncatemacro{\mminusone}{\m-1}

\foreach \k in {0,...,\mminusone}{
  \pgfmathparse{int(2*(\k-\m))}\edef\y{\pgfmathresult}
  \pgfmathparse{int(\k+1)}\edef\kp{\pgfmathresult}
  \node[anchor=east, xshift=-8pt, yshift=-1pt, draw=none] at (0_\y) {\footnotesize$t_{i,\kp}$};
}

\node[anchor=east, xshift=-8pt, yshift=-1pt, draw=none] at (0_0) {\footnotesize$\overline{a}$};

\foreach \k in {1,...,\n}{
  \pgfmathparse{int(2*\k)}\edef\y{\pgfmathresult}
  \node[anchor=east, xshift=-8pt, yshift=-1pt, draw=none] at (0_\y) {\footnotesize$t_{j,\k}$};
}

\foreach \k in {1,...,\m}{
  \pgfmathparse{int(2*(\k-\m)-1)}\edef\y{\pgfmathresult}
  \node[anchor=west, xshift=8pt, yshift=-1pt, draw=none] at (2_\y) {\footnotesize$s_{i,\k}$};
}

\foreach \k in {1,...,\n}{
  \pgfmathparse{int(2*\k-1)}\edef\y{\pgfmathresult}
  \node[anchor=west, xshift=8pt, yshift=-1pt, draw=none] at (2_\y) {\footnotesize$s_{j,\k}$};
}

\foreach \i in {1,...,\n}{
  \pgfmathparse{int(2*\i)}\edef\y{\pgfmathresult};
  \pgfmathparse{int(2*\i-1)}\edef\t{\pgfmathresult};
  \draw [<-] (0_\y.-10) to (2_\t.140);
  \draw [<-] (0_\y.-35) to (2_\t.165);
}

\foreach \i in {2,...,\n}{
  \pgfmathparse{int(2*\i)}\edef\y{\pgfmathresult};
  \pgfmathparse{int(2*\i-3)}\edef\t{\pgfmathresult};
  \draw [->] (0_\y) to (2_\t);
}

\foreach \i in {1,...,\n}{
  \pgfmathparse{int(2*\i-1)}\edef\t{\pgfmathresult};
  \pgfmathparse{int(2*\i-2)}\edef\y{\pgfmathresult};
  \draw [->] (0_\y) to (2_\t);
}

\draw [<-] (0_0) to (0_2);
\draw [<-] (0_0) to (0_-2);

\foreach \i in {1,...,\m}{
  \pgfmathparse{int(-2*\i)}\edef\y{\pgfmathresult};
  \pgfmathparse{int(1-2*\i)}\edef\t{\pgfmathresult};
  \draw [<-] (0_\y.10) to (2_\t.-140);
  \draw [<-] (0_\y.35) to (2_\t.-165);
}

\foreach \i in {2,...,\m}{
  \pgfmathparse{int(-2*\i)}\edef\y{\pgfmathresult};
  \pgfmathparse{int(3-2*\i)}\edef\t{\pgfmathresult};
  \draw [->] (0_\y) to (2_\t);
}

\foreach \i in {1,...,\m}{
  \pgfmathparse{int(1-2*\i)}\edef\t{\pgfmathresult};
  \pgfmathparse{int(2-2*\i)}\edef\y{\pgfmathresult};
  \draw [->] (0_\y) to (2_\t);
}
\end{tikzpicture}

\end{minipage}

\caption{Seed $\mu_a(\quiver_a) = \quiver_{\overline a}$ for $\dimv_i=4$, $\dimv_j=3$.}
\label{fig-mut-Qmn}
\end{minipage}%
\end{figure}

Next, we use the same set of integers $\{1,\ldots, 2(\dimv_i+\dimv_j)-3\}$ to populate a $(2\dimv_i-1)\times (2\dimv_j-1)$ rectangular array as shown in Figure~\ref{fig-array}: we label the southern corner of the array by~1, with the labels increasing by 1 for each step north-east and north-west. The sequence $\bs \mu_{ a}$ has one mutation for each entry in this array. To describe the order in which these mutations are performed, we read the array north-south column by north-south column starting from the west, where in each such column we first read the circled entries from south to north, and then read the uncircled entries from north to south. For example, reading the array in Figure~\ref{fig-array} this way produces the following sequence of mutations starting from the seed $\quiver_{a:i\rightarrow j}$ for $\dimv_i=4$, $\dimv_j=3$:
\beq
\label{eq:Baxter-seq}
\bs\mu_{a} = \bs\mu_{a}^{(11)} \circ \ldots \circ \bs\mu_{a}^{(1)},
\eeq
where
\begin{align*}
&\bs\mu_{a}^{(1)} = \mu_5, &
&\bs\mu_{a}^{(2)} = \mu_6 \circ \mu_4, &
&\bs\mu_{a}^{(3)} = \mu_5 \circ \mu_7\mu_3, \\
&\bs\mu_{a}^{(4)} = \mu_4\mu_8 \circ \mu_6\mu_2, &
&\bs\mu_{a}^{(5)} = \mu_3\mu_7 \circ \mu_9\mu_5\mu_1, &
&\bs\mu_{a}^{(6)} = \mu_2\mu_6\mu_{10} \circ \mu_8\mu_4, \\
&\bs\mu_{a}^{(7)} = \mu_5\mu_9 \circ \mu_{11}\mu_7\mu_3, &
&\bs\mu_{a}^{(8)} = \mu_4\mu_8 \circ \mu_{10}\mu_6, &
&\bs\mu_{a}^{(9)} = \mu_7 \circ \mu_9\mu_5, \\
&\bs\mu_{a}^{(10)} = \mu_6 \circ \mu_8, &
&\bs\mu_{a}^{(11)} = \mu_7.
\end{align*}

\begin{figure}[h]
\begin{tikzpicture}[every node/.style={inner sep=0, minimum size=0.5cm, circle}, x=0.7cm, y=0.35cm]

\def\m{3};
\def\n{2};

\pgfmathparse{int(2*\m)};
\edef\a{\pgfmathresult};
\pgfmathparse{int(2*\n)};
\edef\b{\pgfmathresult};
\foreach \i in {0,...,\a}
	\foreach \j in {0,...,\b}
	{
		\pgfmathparse{int(\i+\j)};
		\edef\x{\pgfmathresult};
		\pgfmathparse{int(\i-\j)};
		\edef\y{\pgfmathresult};
		\pgfmathparse{int(2*\n+\y+1)};
		\edef\v{\pgfmathresult};
	
		\pgfmathparse{int(mod(\i,2))};
		\edef\t{\pgfmathresult};
		\ifnum \t=0
			\node[draw] (\x_\y) at (\x,\y) {\tiny{$\v$}};
		\else
			\node (\x_\y) at (\x,\y) {\tiny{$\v$}};
		\fi

		\pgfmathparse{int(mod(\i+\j,2))};
		\edef\s{\pgfmathresult};
		\ifnum \s=1
			\draw[dashed, thin] (\i+\j-.5,-\b-1) to (\i+\j.5,-\b-1) to (\i+\j.5,\a+1) to (\i+\j-.5,\a+1) to (\i+\j-.5,-\b-1);
		\fi
	}

\end{tikzpicture}

\caption{Mutation sequence $\bs \mu_{ a}$ for $\dimv_i=4$, $\dimv_j=3$.}
\label{fig-array}
\end{figure}

The mutation sequences $\bs \mu_{a}$ were studied in~\cite{SS25} in the context of special coordinate systems in higher Teichm\"uller theory. In particular, we recall the following fact from that paper:

\begin{lemma}
\label{lem:Baxter-trop-var}
Each sign~\eqref{basis.mutatation} in the mutation sequence $\bs \mu_{ a}$ is positive. The resulting seed $\mu_a(\quiver_a)=\{e'_m\}$ is given by 
$$
e'_m =
\begin{cases}
e_{m+2\dimv_j-1} &\text{if} \;\; m <2\dimv_i-1, \\
-\sum_{r=1}^{2(\dimv_i+\dimv_j)-3}e_r  &\text{if} \;\; m = 2\dimv_i-1, \\
e_{m-2\dimv_i+1} &\text{if} \;\; m >2\dimv_i-1.
\end{cases}
$$
\end{lemma}

Since $e_{\overline a}=e'_{2\dimv_i-1}$, we get the following result.
\begin{cor}
\label{cor:arrow-reverse}
We have
$$
\Pi_{\Gamma_{\overline a}}= \bs \mu_{ a}(\Pi_{\Gamma}).
$$
\end{cor}

\begin{defn}
\label{def:bibax}
If $ a$ is an arrow between distinct vertices of $\Gamma$, we define the corresponding \emph{bi-fundamental Baxter operator} to be the inverse of the product of quantum dilogarithm power series
$$
\Qop_ a = \Psi_{\quiver';\quiver} \in \widehat{\mathcal{R}}_\quiver
$$
associated to $\bs \mu_{ a}$ by the formula~\eqref{eq:qdls}, where $\quiver=\quiver_\Gamma$ and $\quiver'=\quiver_{\Gamma_{\overline a}}$.
\end{defn}

We recall from Theorem~\ref{thm:keller} that the formal series $\Qop_a$ is completely determined by the pair of seeds $\quiver$, $\quiver'$ and is independent of the  combinatorial factorization of the cluster transformation $\bs\mu_a$ we presented earlier. As a direct corollary of Lemma~\ref{lem:Baxter-trop-var}, we get:
\begin{cor} 
\label{frozen-evolution}
Let $a \colon i\rightarrow j$ be an arrow between distinct vertices of $\Gamma$. Writing $\Qop_a$ for the corresponding bi-fundamental Baxter operator, we have
$$
\Lbb_{\quiver_{\Gamma_{\overline a}}} = \Ad_{\Qop_{a}} \left(\Lbb_{\quiver_\Gamma}\right) \subset \Lbb_{G_V,T_F,0}.
$$
\end{cor}

Another useful fact about the power series $\mathbf{Q}_a$ we recall from~\cite{SS25} is that directions $f_r$ of the monomial quantum dilogarithm arguments $Y_{f_r}$ in the factorization
\beq
\label{eq:Baxter-factor}
\mathbf{Q}_ a = \Psi\big(Y_{f_{(2\dimv_i-1)(2\dimv_j-1)}}\big) \cdots \Psi\big(Y_{f_1}\big)
\eeq
satisfy
\beq
\label{eq:positivity}
\big(f_r,x_{\varpi_{j,\dimv_j}}\big) >0, \quad \big(f_r,x_{\varpi_{i,\dimv_i}}\big)<0, \qquad 1\leq r \leq (2\dimv_i-1)(2\dimv_j-1).
\eeq

Now let us introduce an additional formal parameter $t$ and consider the  space $\Fcr^\Gamma_\Q\brr{t}$ of $K_{T_F}(\pt)\otimes_\Z \Q(q)\brr{t}$-valued functions $f$ on the weight lattice $\Prm$ of $G_V$ such that for each integer $k$ the function $f_k$ obtained by taking the coefficient of $t^k$ in $f$ has finite support. For a node $j \in \Gamma_0^\rg$ we define a representation $\varrho^{(j)}_t$ of the quantum torus $\Tbb_{\Lambda_{|\Gamma|}}$ on $\Fcr^\Gamma_\Q\brr{t}$ by twisting the standard representation~\eqref{eq:toda-torus-F-action} by the automorphism of $\Tbb_{\Lambda_{|\Gamma|}}$ which is the identity on all generators except for the momenta $P_{j,n}$, which we rescale by $P_{j,n}\mapsto tP_{j,n}$ for all $1\leq n\leq \dimv_j$, so that
$$
(P_{j,n}f)(\la) = tf(\la-\eps_{j,n}).
$$
In the case $j \in \Gamma_0^\rf$ is a framing node, we define the representation $\varrho^{(j)}_t$ in a similar way by rescaling the generator $z_j$ by $t$, so that
$$
(z_jf)(\la) = tz_jf(\la).
$$
Then it follows from~\eqref{eq:positivity} that for each arrow $a \colon i \to j$ the series $\Qop_ a$ determines a well-defined operator $\Qop_ a(t)=\varrho_t^{(j)}(\Qop_ a)$ on $\Fcr^\Gamma_\Q\brr{t}$. We also note that the representation $\varrho_1 = \varrho_1^{(j)}$ given by the non-rescaled formulas~\eqref{eq:toda-torus-F-action} extends to a representation of $\Tbb_{\Lambda_{|\Gamma|}}\bss{t}$ on $\Fcr^\Gamma_\Q\brr{t}$.

Similarly, we define the $t$-completed space of distributions $\Dcr^\Gamma_\Q\brr{t}$, and write
$$
\Dcr^\Gamma_{\Q;\bs w}\brr{t} = \Dcr^\Gamma_\Q\brr{t} \otimes_{\Q(q)} \Kc^\Gamma_{\Q;\bs w}, \qquad \Kc^\Gamma_{\Q;\bs w} = \Kc^\Gamma_{\bs w} \otimes_{\Z[q^{\pm1}]} \Q(q).
$$
Setting $\bar t = t$, we extend the pairing~\eqref{eq:wwpair} to
\beq
\label{eq:wwpair-t}
(\cdot, \cdot) \colon \Fcr^\Gamma_\Q\brr{t} \otimes_{\Q(q)\brr{t}} \Dcr^\Gamma_{\Q;\bs w}\brr{t} \longrightarrow \Kc^\Gamma_{\Q;\bs w}.
\eeq
We also have a Whittaker transform isomorphism
\beq
\label{eq:alg-whit-gamma-t}
\Wc_t^\Gamma \colon \Vcr^{\Gamma}_\Q\brr{t} \simeq  K_{G_V\times T_F}(\pt)\otimes_{\mathbb{Z}}\mathbb{Q}(q)\brr{t},
\eeq
obtained from~\eqref{eq:alg-whit-gamma} by extending scalars to $\mathbb{Q}(q)$ and completing in $t$.

Now we recall the following result from~\cite{SS25}.

\begin{prop}[\cite{SS25} Proposition 7.10]
\label{prop:bax-eig}
Let $\Qop_a^*(t)$ be the formal adjoint of $\Qop_ a(t)$ with respect to the pairing~\eqref{eq:wwpair-t}. Then the $G_V$-Whittaker kernel diagonalizes the operator $\Qop_a^*(t)$:
\beq
\label{eq:bax-eig}
\Qop_a^*(t) \circ W^{\Gamma} = \mathbf{\Psi}_a(t)^{-1} \cdot W^\Gamma,
\eeq
where $\mathbf{\Psi}_a(t)$ is the $t$-twisted version of formula~\eqref{eq:Psi_a}:
$$
\mathbf{\Psi}_a(t)= 
\begin{cases}
\prod_{r=1}^{\dimv_i}\prod_{s=1}^{\dimv_j}\Psi\left(tu(a) w_{i,r}/w_{j,s}\right) &\text{if} \; a \in \Gamma_1^\rg, \\
\prod_{s=1}^{\dimv_j}\Psi\left(tz_i/w_{j,s}\right) &\text{if} \; a \in \Gamma_1^{\rf_+}, \\
\prod_{r=1}^{\dimv_i}\Psi\left(tw_{i,r}/z_j\right) &\text{if} \; a \in \Gamma_1^{\rf_-}.
\end{cases}
$$
\end{prop}

Note that the product on the right-hand-side of~\eqref{eq:bax-eig} is a power series in $t$ 
whose coefficients are $S_{\dimv_i}\times S_{\dimv_j}$ symmetric Laurent polynomials in variables $\bs w_{i,\bullet}, \bs w_{j,\bullet}$ with coefficients in $\mathbb{Q}(q)[u_a^{\pm1}, z_k^{\pm1}]$. On the other hand, the Pieri rule~\eqref{eq:chi-eig} tells us that the action of the Toda Hamiltonians on the Whittaker kernel corresponds to multiplication by elementary symmetric polynomials. Since the elementary symmetric polynomials generate the ring of all symmetric polynomials, it follows that there is an element
$$
A_a(t) \in \Lbb_{GL_{\dimv_i}\times GL_{\dimv_j},0}\bss{t}\otimes_{\mathbb{Z}[q^{\pm1}]}\mathbb{Q}(q)
$$
whose action via the un-rescaled representation $\varrho_1$ on $W^\Gamma$ produces the same eigenvalue.  Since no nonzero element of the quantum torus can annihilate the Whittaker kernel, it follows that $\Qop_ a(t)=\varrho_1(A_a(t))$ and so by Lemma~\ref{lem:pres} we deduce that $\Qop_ a(t)$ preserves the subspace $\Vcr^\Gamma_\Q\brr{t} \subset \Fcr^\Gamma_\Q\brr{t}$ of functions vanishing outside the dominant cone.

\begin{prop}
\label{prop:inter}
Let $ a\colon i\rightarrow j$ be an arrow between distinct vertices of $\Gamma$. Then for all $A\in \Lbb_{\Pi_\Gamma}$, we have 
\begin{align}
\label{eq:intereq}
\mathbb{W}_\Gamma(\Ad_{\Qop_ a}(A)) = \Ad_{\mathbf{\Psi}_a}\left(\mathbb{W}_\Gamma(A)\right).
\end{align}
In particular, $\Ad_{\mathbf{\Psi}_a}\left(\mathbb{W}_\Gamma(A)\right)$ lies in $\Acr_{G_V,T_F,0}$.
\end{prop}
\begin{proof}
For $j \in \Gamma_0^\rg$, consider the faithful $q$-difference operator representation $\eta^{(j)}_t$ of $\Acr_{G_V,T_F,0}$ on the space $K_{G_V\times T_F}(\pt)\otimes_{\mathbb{Z}}\mathbb{Q}(q)\brr{t}$ defined by precomposing the standard action with the automorphism which rescales $w_{j,n}\mapsto tw_{j,n}$ for all $1\leq n\leq v_j$ and is the identity on all other generators. If $j \in \Gamma_0^\rf$ is the framing node, the representation $\eta^{(j)}_t$ is defined similarly, except in that case we rescale $z_j\mapsto tz_j$. The twisted action $\rho^{(j)}_t$ gives a faithful representation of the algebra $\Lbb_{G_V,T_F,0}$ on $\Vcr^{\Gamma}_\Q\brr{t}$, and as in the proof of Proposition~\ref{prop:alg-whit} the isomorphism $\mathcal{W}_t^\Gamma$ defined by~\eqref{eq:alg-whit-gamma-t} identifies $\Lbb_{G_V,T_F,0}$ and $\Acr_{G_V,T_F,0}$ as subrings in its endomorphism algebra. Moreover, by comparing images of the generators for $\Lbb_{G_V,T_F,0}$ it is easy to see that this identification, constructed with the help of the $t$-twisted Whittaker transform~\eqref{eq:alg-whit-gamma-t}, coincides with the isomorphism~\eqref{eq:boldW-Gamma} constructed with the help of~\eqref{eq:alg-whit-gamma}, i.e.\ for all $A\in \Lbb_{G_V,T_F,0}$ we have
$$
\mathcal{W}_t^\Gamma\circ\varrho_t(A) = \eta_t(\mathbb{W}_\Gamma(A))\circ\mathcal{W}_t^\Gamma.
$$ 
On the other hand, by Proposition~\ref{prop:bax-eig} we have
\begin{align}
\label{eq:qop-inter}
\mathcal{W}_t^\Gamma\circ \varrho_t(\Qop_ a) = \eta_t(\mathbf{\Psi}_{a})\circ \mathcal{W}_t^\Gamma
\end{align}
and so~\eqref{eq:intereq} follows. 
\end{proof}
Note that as a corollary of the relation~\eqref{eq:qop-inter} intertwining the bi-fundamental Baxter operator $\Qop_ a$ with a multiplication operator, we obtain

\begin{cor}
\label{cor:Baxter-comm}
The bi-fundamental Baxter operators $\Qop_{ a_1},\Qop_{ a_2}$ associated to any pair of arrows $ a_1, a_2$ commute.
\end{cor}

\subsection{Isomorphism with the Coulomb branch algebra}

Our goal in this section is to prove the following result:
\begin{theorem}
\label{thm:main-iso}
Let $\Gamma$ be a quiver without 1-cycles and $\quiver_\Gamma\subset\Lambda_{|\Gamma|}$ the cluster seed constructed in Definition~\ref{eq:seed-def}. Then $G_V$-Whittaker transform $\mathbb{W}_{\Gamma}\colon\Lbb_{G_V,T_F,0} \simeq \Acr_{G_V,T_F,0}$ restricts to an algebra isomorphism
\beq
\label{eq:main-iso}
\mathbb{W}_\Gamma\colon \Lbb_\Gamma \simeq  \Acr_\Gamma,
\eeq
where we abbreviate $\Lbb_\Gamma = \Lbb_{\quiver_\Gamma}$.
\end{theorem}

Here is the scheme of the proof of Theorem~\ref{thm:main-iso}. Given an arrow $ a\in \Gamma_1$,  write 
$$
\quiver_{G_V,T_F,N_a}=\quiver_{G_V,T_F,0}\sqcup\{e_ a\}
$$ 
with vector $e_a$ defined by~\eqref{eq:earr}. Then by the characterization~\eqref{quantum.cluster.ub} of $\Lbb_\Pi$, we have
$$
\Lbb_\Gamma = \bigcap_{a\in \Gamma_1} \Lbb_{\Gamma,a},
$$
where $\Lbb_{\Gamma,a} = \Lbb_{\quiver_{G_V,T_F,N_a}} \subset \Lbb_{G_V,T_F,0}$ is the subalgebra consisting of elements that remain Laurent under mutation in direction $e_ a$. Recalling Proposition~\ref{prop:conj-crit} and Remark~\ref{rmk:coprime} we also have
$$
\Acr_\Gamma = \bigcap_{ a\in \Gamma_1}  \Acr_{\Gamma, a},
$$
where
$$
\Acr_{\Gamma, a}=\left\{A\in \Acr_{G_V,T_F,0} ~\big|~ \Ad_{\mathbf{\Psi}_a}(A) \in \Acr_{G_V,T_F,0} \right\}.
$$

To prove Theorem~\ref{thm:main-iso} it suffices to show that for each arrow $a \in \Gamma_1$ the $G_V$-Whittaker transform $\mathbb{W}_{\Gamma}\colon\Lbb_{G_V,T_F,0} \simeq \Acr_{G_V,T_F,0}$ restricts to an isomorphism
\begin{align}
\label{eq:intersection-iso}
\mathbb{W}_\Gamma \colon \Lbb_{\Gamma,a} \simeq \Acr_{\Gamma, a}.
\end{align}
Given an arrow $a \colon i\rightarrow j$ note that the element $Y_{e_ a} \in \Tbb_{\La_{|\Gamma|}}$ commutes with all $P_{k,r},X_{k,r} \in \Tbb_{\Lambda_{|\Gamma|}}$ unless $k\in \{i,j\}$, and similarly $\mathbf{\Psi}_a$ commutes with all $D_{k,r}, w_{k,r} \in \Acr_{T_V,T_F,0}$ unless $k\in \{i,j\}$. So using the splitting~\eqref{eq:seed-split}, we see that ~\eqref{eq:intersection-iso} will follow in general provided we can prove Theorem~\ref{thm:main-iso} in the special case that $\Gamma=i\rightarrow j$ consists of a pair of nodes $i,j$ and a single arrow $a \colon i\rightarrow j$. {We now do this in the case that both $i,j$ are gauge nodes, i.e.\ $a \in \Gamma_1^\rg$. The case $a \in \Gamma_1^\rf$, is treated by a simpler version of the argument below.}

Recall the seed $\quiver_{i\rightarrow j} = \quiver_a$ in the lattice $\Lambda_{i-j}=\Lambda_{GL_{\dimv_i}}\oplus \Lambda_{GL_{\dimv_j}}$, defined in~\eqref{eq:arrow-subseeds}. Our goal is to show that
\begin{align}
\label{eq:ij-iso}
\mathbb{W}_{i \to j} \colon \Lbb_{i \to j} \simeq \Acr_{i \to j}.
\end{align}
We will do this by constructing an intermediate algebra $\Lbb^{\mathrm{res}}_{i \to j^\diamondsuit} \subset \Lbb_{GL_{\dimv_i},0} \otimes_{\Z[q^{\pm1}]} \Acr_{GL_{\dimv_j},0}$, and showing that the algebra isomorphisms
\begin{align*}
\Wbb_{(j)} \colon \Lbb_{GL_{\dimv_i},0} \otimes_{\mathbb{Z}[q^{\pm1}]} \Lbb_{GL_{\dimv_j},0} &\simeq \Lbb_{GL_{\dimv_i},0} \otimes_{\mathbb{Z}[q^{\pm1}]} \Acr_{GL_{\dimv_j},0}, \\
\Wbb_{(i)} \colon \Lbb_{GL_{\dimv_i},0} \otimes_{\mathbb{Z}[q^{\pm1}]} \Acr_{GL_{\dimv_j},0} &\simeq \Acr_{GL_{\dimv_i},0} \otimes_{\mathbb{Z}[q^{\pm1}]} \Acr_{GL_{\dimv_j},0}
\end{align*}
restrict to isomorphisms
$$
\Wbb_{(j)} \colon \Lbb_{i \to j} \simeq \Lbb^{\mathrm{res}}_{i \to j^\diamondsuit}, \qquad \Wbb_{(i)} \colon \Lbb^{\mathrm{res}}_{i \to j^\diamondsuit} \simeq \Acr_{i \to j}.
$$
The first step is to consider the cluster seed $\quiver_{i \to j^\diamondsuit}$ associated to the quiver $i \to j^\diamondsuit$ obtained from $i \to j $ by replacing the gauge node $j$ with a collection of $\dimv_j$ gauge nodes, each carrying a 1-dimensional vector space. In other words, we consider the theory with the same matter content, but in which the $GL_{\dimv_j}$ gauge symmetry is broken down to its maximal torus $T_{\dimv_j}\subset GL_{\dimv_j}$, {as in the discussion in Remark~\ref{rem:G-to-T-sym}}. More precisely, we form the lattice
$$
\Lambda_{i - j^\diamondsuit} = \Lambda_{GL_{\dimv_i}}\oplus \bigoplus_{m=1}^{\dimv_j}\mathbb{Z}\langle \beta_{j,m}, \gamma_{j,m}\rangle
$$
with the skew-form making the direct summands orthogonal and such that
$$
(\gamma_{j,m},\beta_{j,n})=\delta_{m,n}, \qquad (\beta_{j,m},\beta_{j,n})=0=(\gamma_{j,m},\gamma_{j,n})
$$
for all $1 \le m,n \le \dimv_j$. Sending
$$
Y_{\beta_{j,s}} \mapsto D_{j,s}, \qquad Y_{\gamma_{j,s}} \mapsto w_{j,s}
$$
we identify the quantum torus $\Tbb_{i \to j^\diamondsuit}$ with $\Tbb_{GL_{\dimv_i},0} \otimes \Dc(T_{\dimv_j})$, where 
$$
\Dc(T_{\dimv_j}) = \Z[q^{\pm1}]\ha{D_{j,r},w_{j,r} \,|\, 1 \le r \le d}/\langle D_rw_s = q^{2\delta_{r,s}}w_sD_r\rangle
$$
is the algebra of $q$-difference operators on the torus $T_{\dimv_j}$ from Section~\ref{subsec:localization}. We equip the lattice $\Lambda_{i - j^\diamondsuit}$ with the seed
\begin{align}
\label{eq:seedbase}
\quiver_{i\rightarrow j^\diamondsuit} = \{e_{s_{i,m}},e_{t_{i,m}}\}_{m=1}^{v_i-1} \sqcup \{e_{a_n}=\gamma_{j,n}-p_{i,1}\}_{n=1}^{\dimv_j}
\end{align}
which we represent graphically in Figure~\ref{fig:Qab}. 

\begin{figure}[h]

\begin{tikzpicture}[thick, x=.7cm, y=1cm, every node/.style={circle, draw, fill=white, minimum size=6mm, inner sep=0pt, font=\footnotesize}]

\node (s1) at (-1,-1) {$s_{i,1}$};
\node (t1) at (-1,1) {$t_{i,1}$};
\node (s2) at (-3,-1) {$s_{i,2}$};
\node (t2) at (-3,1) {$t_{i,2}$};
\node (s3) at (-5,-1) {$s_{i,3}$};
\node (t3) at (-5,1) {$t_{i,3}$};

\node (v1) at (1,0) {$a_1$};
\node (v2) at (3,0) {$a_2$};
\node (v3) at (5,0) {$a_3$};

\draw[->] (v1) to[out=-120, in=30] (s1);
\draw[->] (v2) to[out=-135, in=15] (s1);
\draw[->] (v3) to[out=-150, in=0] (s1);

\draw[->] (t1) to (s2);
\draw[->] (t2) to (s3);

\draw[->] (t3) to (s2);
\draw[->] (t2) to (s1);

\draw[->] (t1) to[out=-30, in=120] (v1);
\draw[->] (t1) to[out=-15, in=135] (v2);
\draw[->] (t1) to[out=0, in=150] (v3);

\draw[->] (s1.80) to (t1.-80);
\draw[->] (s1.100) to (t1.-100);
\draw[->] (s2.80)  to (t2.-80);
\draw[->] (s2.100) to (t2.-100);
\draw[->] (s3.80)  to (t3.-80);
\draw[->] (s3.100) to (t3.-100);

\end{tikzpicture}

\caption{The seed $\quiver_{i\rightarrow j^\diamondsuit}$ for $\dimv_i=4$, $\dimv_j=3$.}
\label{fig:Qab}
\end{figure}

Now consider the localization $\Tbb^{\loc}_{i\rightarrow j^\diamondsuit}$ of the quantum torus $\Tbb_{i\rightarrow j^\diamondsuit}$ at the Ore denominator set
$$
Ø(j) = \hc{1-q^{2k}w_{j,\alpha} ~\Big|~ \alpha\in\Delta_+(GL_{\dimv_j}), \, k\in\mathbb{Z}}.
$$
We have an algebra isomorphism
$$
\Tbb^\loc_{i\rightarrow j^\diamondsuit}\simeq \Tbb_{GL_{\dimv_i},0} \otimes_{\mathbb{Z}[q^{\pm1}]} \Dc_{\loc}(T_{\dimv_j}),
$$
where $\Dc_{\loc}(T_{\dimv_j})$ is the localization of $\Dc(T_{\dimv_j})$ at $Ø(j)$. Recall Definition~\ref{def:Dres} of the residue subalgebra $\Dres(T_{\dimv_j}) \subset \Dc_\loc(T_{\dimv_j})$ along with isomorphism $\Dres(T_{\dimv_j}) \simeq \Acr_{GL_{\dimv_j},0}$ from Proposition~\ref{rmk:minuscule-res}. Then we have the embedding
$$
\Lbb_{GL_{\dimv_i}}\otimes_{\mathbb{Z}[q^{\pm1}]} \Acr_{GL_{\dimv_j,0}} \subset \Tbb^{\loc}_{i\rightarrow j^\diamondsuit}.
$$
Observe that the set $Ø(j)$ is mutation-invariant, since from~\eqref{eq:seedbase} it is clear that its elements commute with all mutable generators of $\Tbb_{i \to j^\diamondsuit}$.

\begin{defn}
\label{def:Lres}
We define the \emph{residue universally Laurent algebra} $\Lbb^{\res}_{i\rightarrow j^\diamondsuit}$ to be the subalgebra of all elements 
$$
{B \in \Lbb_{GL_{\dimv_i}}\otimes_{\mathbb{Z}[q^{\pm1}]} \Acr_{GL_{\dimv_j,0}}}
$$
such that  $fB\in \Lbb_{i\rightarrow j^\diamondsuit}$ for some {$f\in K_{GL_{\dimv_j} \times \widetilde{\mathbb{C}}^\times}(\pt)$}.
\end{defn}

\begin{remark}
\label{rmk:laur-crit}
Using the characterization~\eqref{quantum.cluster.ub} and the fact that elements of $Ø(j)$ commute with mutations, it follows that an element $B \in \Lbb_{GL_{\dimv_i},0}\otimes \Acr_{GL_{\dimv_j},0}$ lies in the residue universally Laurent algebra $\Lbb^{\res}_{i\rightarrow j^\diamondsuit}$ if and only if for all $1 \le n \le \dimv_j$ we have $\mu_{a_n}(B) \in \Tbb^{\loc}_{i\rightarrow j^\diamondsuit}$.
\end{remark}

\begin{prop}
\label{prop:whit-j}
The algebra isomorphism
$$
\Wbb_{(j)} \colon \Lbb_{GL_{\dimv_i},0} \otimes_{\mathbb{Z}[q^{\pm1}]} \Lbb_{GL_{\dimv_j},0} \simeq \Lbb_{GL_{\dimv_i},0} \otimes_{\mathbb{Z}[q^{\pm1}]} \Acr_{GL_{\dimv_j},0}
$$
restricts to an isomorphism
$$
\mathbb{W}_{(j)} \colon \Lbb_{i\rightarrow j} \simeq \Lbb^{\res}_{i\rightarrow j^\diamondsuit}.
$$
\end{prop}

\begin{proof}
Let us define a sequence of mutations $\bs\mu_{\Box \to j}$ in $\quiver_{i \to j}$, which in the numbering convention of Section~\ref{subsec:bibax} reads
$$
\bs\mu_{\Box \to j} = \mu_1 \circ \mu_2 \circ \ldots \circ \mu_{2\dimv_j-1}.
$$
We denote by $\Qop_{\Box \to j} \in \widehat\Rc_{i \to j}$ the corresponding Baxter operator. Note that if a quantum torus element remains Laurent under a sequence of mutations in distinct directions, it must in particular remain Laurent under the first mutation in the sequence. Since by~\eqref{quantum.cluster.ub} we know that the Laurent algebra $\Lbb_{i \to j}$ is the subalgebra of those elements in $\Lbb_{GL_{\dimv_i},0} \otimes_{\Z[q^{\pm1}]} \Lbb_{GL_{\dimv_j},0}$ that remain Laurent under the mutation in direction $e_a$, the previous observation yields the following equivalent description:
$$
\Lbb_{i \to j} = \hc{A \in \Lbb_{GL_{\dimv_i},0} \otimes_{\Z[q^{\pm1}]} \Lbb_{GL_{\dimv_j},0} ~\Big|\, \Ad_{\Qop_{\Box \to j}}(A) \in \Tbb_{i \to j}}.
$$
On the other hand, by Lemma~\ref{lem:Baxter-trop-var} applied in the case $\dimv_i=1$, the vectors $\hc{e_{s_{j,n}},e_{t_{j,n}}}$ belong to $\bs\mu_{\Box \to j}(\quiver_{i \to j})$ for all $1 \le n \le d_j$. Therefore for any $A \in \Lbb_{i \to j}$, the element $\Ad_{\Qop_{\Box \to j}}(A)$ lies in $\Tbb_{GL_{\dimv_i},0} \otimes_{\Z[q^{\pm1}]} \Lbb_{GL_{\dimv_j},0} \subset \Tbb_{i \to j}$, so we can further reformulate the description of $\Lbb_{i \to j}$ above as
$$
\Lbb_{i \to j} = \hc{A \in \Lbb_{GL_{\dimv_i},0} \otimes_{\Z[q^{\pm1}]} \Lbb_{GL_{\dimv_j},0} ~\Big|\, \Ad_{\Qop_{\Box \to j}}(A) \in \Tbb_{GL_{\dimv_i},0} \otimes_{\Z[q^{\pm1}]} \Lbb_{GL_{\dimv_j},0}}.
$$
Next, consider the sequence $\bs\mu_{\Box \to j^\diamondsuit} = \mu_{a_1} \circ \ldots \circ \mu_{a_{\dimv_j}}$ of $\dimv_j$ commuting mutations in $\quiver_{i\rightarrow j^\diamondsuit}$, and the corresponding product of quantum dilogarithms
$$
\bs\Psi_{\Box \to j^\diamondsuit} = \prod_{r=1}^{v_j}\Psi(Y_{p_{i,1}-\gamma_{j,r}}).
$$
Note that a quantum torus element remains Laurent under a sequence of mutations in distinct, commuting directions if and only if it remains Laurent under mutation in each direction separately. Thus by Remark~\ref{rmk:laur-crit} we have the following description of the algebra $\Lbb^{\res}_{i \to j^\diamondsuit}$:
$$
\Lbb^{\res}_{i \to j^\diamondsuit} = \hc{B \in \Lbb_{GL_{\dimv_i},0} \otimes_{\Z[q^{\pm1}]} \Acr_{GL_{\dimv_j},0} ~\Big|\, \Ad_{\bs\Psi_{\Box \to j^\diamondsuit}}(B) \in \Tbb^{\loc}_{i \to j^\diamondsuit}}.
$$
On the other hand, we observe that for all $B\in \Lbb^{\res}_{i \to j^\diamondsuit}$, the element $\Ad_{\bs\Psi_{\Box \to j^\diamondsuit}}(B)$ belongs to $\Tbb_{GL_{\dimv_i,0}} \otimes_{\Z[q^{\pm1}]} \Dc_{\res}(T_{\dimv_j}) \subset \Tbb^\loc_{i \to j^\diamondsuit}$. Indeed, the $S_{\dimv_j}$ symmetry of $\Ad_{\bs\Psi_{\Box \to j^\diamondsuit}}(B)$ is immediate from that of $B$ and $\bs\Psi_{\Box \to j^\diamondsuit}$, and that it satisfies the residue conditions~\eqref{eq:residue-condition} follows from Lemma~\ref{rmk:res-pres}. Therefore, we obtain
$$
\Lbb^{\res}_{i \to j^\diamondsuit} = \hc{B \in \Lbb_{GL_{\dimv_i},0} \otimes_{\Z[q^{\pm1}]} \Acr_{GL_{\dimv_j},0} ~\Big|\, \Ad_{\bs\Psi_{\Box \to j^\diamondsuit}}(B) \in \Tbb_{GL_{\dimv_i,0}} \otimes_{\Z[q^{\pm1}]} \Acr_{GL_{\dimv_j},0}}.
$$
Finally, by Proposition~\ref{prop:inter} we have
$$
\Wbb_{(j)}(\Ad_{\Qop_{\Box \to j}}(A)) = \Ad_{\bs\Psi_{\Box \to j^\diamondsuit}}(\Wbb_{(j)}(A)).
$$
Setting $B = \Wbb_{(j)}(A)$, we see that
$$
\Ad_{\Qop_{\Box \to j}}(A) \in \Tbb_{GL_{\dimv_i},0} \otimes_{\Z[q^{\pm1}]} \Lbb_{GL_{\dimv_j},0} \iff \Ad_{\bs\Psi_{\Box \to j^\diamondsuit}}(B) \in \Tbb_{GL_{\dimv_i,0}} \otimes_{\Z[q^{\pm1}]} \Acr_{GL_{\dimv_j},0},
$$
hence
$$
\Wbb_{(j)} \colon \Lbb_{i \to j} \simeq \Lbb^{\res}_{i \to j^\diamondsuit}.
$$
\end{proof}

\begin{prop}
\label{prop:whit-i}
The algebra isomorphism
$$
\Wbb_{(i)} \colon \Lbb_{GL_{\dimv_i},0} \otimes_{\mathbb{Z}[q^{\pm1}]} \Acr_{GL_{\dimv_j},0} \simeq \Acr_{GL_{\dimv_i},0} \otimes_{\mathbb{Z}[q^{\pm1}]} \Acr_{GL_{\dimv_j},0}
$$
restricts to an isomorphism
$$
\Wbb_{(i)} \colon \Lbb^{\res}_{i\rightarrow j^\diamondsuit} \simeq \Acr_{i \to j}.
$$
\end{prop}
\begin{proof}
We start by giving a description of the residue universal Laurent algebra $\Lbb^{\res}_{i \to j^\diamondsuit}$ similar to the one in Proposition~\ref{prop:whit-j}. Consider the sequence of mutations
$$
\bs\mu_{i \to j^\diamondsuit_n} = \mu_{s_{i,\dimv_i-1}} \circ \mu_{t_{i,\dimv_i-1}} \circ \ldots \circ \mu_{s_{i,1}} \circ \mu_{t_{i,1}} \circ \mu_{a_n},
$$
and let $\Qop_{i \to j^\diamondsuit_n}$ be the corresponding Baxter operator. Again by Lemma~\ref{lem:Baxter-trop-var}, the vectors $\hc{e_{s_{i,m}},e_{t_{i,m}}}$ belong to the seed obtained by applying $\bs\mu_{i \to j^\diamondsuit_n}$, and hence for any $A \in \Lbb^{\res}_{i \to j^\diamondsuit}$, the element $\Ad_{\Qop_{i \to j^\diamondsuit_n}}(A)$ lies in $\Lbb_{GL_{\dimv_i},0} \otimes_{\Z[q^{\pm1}]} \Dc_\loc(T_{\dimv_j}) \subset \Tbb^\loc_{i \to j^\diamondsuit}$. So the residue universal Laurent algebra $\Lbb^{\res}_{i \to j^\diamondsuit}$ consists of all elements $A \in \Lbb_{GL_{\dimv_i},0} \otimes_{\Z[q^{\pm1}]} \Dc_{\res}(T_{\dimv_j})$ such that $\Ad_{\Qop_{i \to j^\diamondsuit_n}}(A) \in \Lbb_{GL_{\dimv_i},0} \otimes_{\Z[q^{\pm1}]} \Dc_\loc(T_{\dimv_j})$ for all $1 \le n \le \dimv_j$.

On the other hand, setting
$$
\bs\Psi_{i \to j^\diamondsuit_n} = \prod_{r=1}^{\dimv_i}\Psi(w_{i,r}/w_{j,n}),
$$
recall that Proposition~\ref{prop:conj-crit} characterizes $\Acr_{i \to j}$ as the subalgebra in $\Dres(T_{\dimv_i}) \otimes_{\Z[q^{\pm1}]} \Dres(T_{\dimv_j})$ consisting of elements $B$ such that $\Ad_{\bs\Psi_{i \to j^\diamondsuit_n}}(B) \in \Dres(T_{\dimv_i}) \otimes_{\Z[q^{\pm1}]} \Dc_\loc(T_{\dimv_j})$ for all $1 \le n \le \dimv_j$. Again by Proposition~\ref{prop:inter}, we have
$$
\Wbb_{(j)}(\Ad_{\Qop_{i \to j^\diamondsuit_n}}(A)) = \Ad_{\bs\Psi_{i \to j^\diamondsuit_n}}(\Wbb_{(j)}(A)),
$$
and for $B = \Wbb_{(j)}(A)$ we get
$$
\Ad_{\Qop_{i \to j^\diamondsuit_n}}(A) \in \Lbb_{GL_{\dimv_i},0} \otimes_{\Z[q^{\pm1}]} \Dc_\loc(T_{\dimv_j}) \iff \Ad_{\bs\Psi_{i \to j^\diamondsuit_n}}(B) \in \Dres(T_{\dimv_i}) \otimes_{\Z[q^{\pm1}]} \Dc_\loc(T_{\dimv_j}),
$$
which yields
$$
\Wbb_{(i)} \colon \Lbb^{\res}_{i \to j^\diamondsuit} \simeq \Acr_{i \to j}.
$$
\end{proof}

Combining Propositions~\ref{prop:whit-j} and~\ref{prop:whit-i}, we deduce that~\eqref{eq:ij-iso} holds whenever $i,j$ are both gauge nodes. The case in which one of $i,j$ is a framing node is handled by a simpler version of the same argument, since there is no need to pass through the intermediate algebra $\Lbb^{\res}_{i\rightarrow j^\diamondsuit}$. This completes the proof of Theorem~\ref{thm:main-iso}.

\section{DT transformation for quivers without 1-cycles}
\label{sec:DT}

In this short section, we show that the seed $\quiver_\Gamma$ associated to a quiver $\Gamma$ without 1-cycles admits a cluster Donaldson--Thomas transformation in the sense of Keller~\cite{Kel17}.

\begin{defn}[\cite{Kel17,GS18}]
A sequence $\bs\mu$ of mutations is called the \emph{cluster Donaldson--Thomas transformation} of the cluster seed $\quiver$ if $\bs\mu(\quiver) = -\quiver$. 
\end{defn}

Given an arrow $a \in \Gamma_1$, recall the seeds $\quiver_\Gamma$, $\quiver_{\Gamma_{\overline{a}}}$ from Section~\ref{subsec:bibax} and the Baxter sequence of mutations $\bs\mu_a$ defined by~\eqref{eq:Baxter-seq}. By Corollary~\ref{cor:arrow-reverse} we have $\bs\mu_a(\quiver_\Gamma) = \quiver_{\Gamma_{\overline{a}}}$. Denote by
$$
\bs\mu_{\Gamma_1} = \prod_{a \in \Gamma_1} \bs\mu_a
$$
the composition of Baxter sequences for all arrows $a \in \Gamma_1$, where we recall that the factors $\bs\mu_a$ commute by Corollary~\ref{cor:Baxter-comm}. Let $\quiver_{\overline\Gamma}$ be the seed obtained from $\quiver_\Gamma$ by replacing each vector $e_{a: i \to j}$ by
\beq
\label{eq:a-bar}
e_{\overline a} = -\sum_{m=1}^{\dimv_i-1}(e_{s_{i,m}} + e_{t_{i,m}}) - e_a - \sum_{n=1}^{\dimv_j-1}(e_{s_{j,n}} + e_{t_{j,n}}).
\eeq
Then we have
$$
\bs\mu_{\Gamma_1}(\quiver_\Gamma) = \quiver_{\overline\Gamma}.
$$

For each gauge node $i \in \Gamma_0^\rg$ and a positive integer $n \in \Z_{\geq0}$ we now define a mutation sequence $\bs\mu_{i;n}$ by setting $\bs\mu_{i;0} = \id$ and requiring
$$
\bs\mu_{i;n+1} = 
\begin{cases}
\bs\mu_{s_i} \circ \bs\mu_{i;n} &\text{if $n$ is even,} \\
\bs\mu_{t_i} \circ \bs\mu_{i;n} &\text{if $n$ is odd,} \\
\end{cases}
$$

where
$$
\bs\mu_{s_i} = \bs\mu_{s_{i,{\dimv_i-1}}} \circ \ldots \circ \bs\mu_{s_{i,1}}, \qquad
\bs\mu_{t_i} = \bs\mu_{t_{i,{\dimv_i-1}}} \circ \ldots \circ \bs\mu_{t_{i,1}}
$$
are the sequences of $\dimv_i-1$ commuting mutations of the sub-seed $\quiver_{GL_{\dimv_i},0} = \hc{e_{s_{i,r}}, e_{t_{i,r}}}_{r=1}^{\dimv_i-1}$. We then form the product
$$
\bs\mu_{\Gamma_0} = \prod_{i \in \Gamma_0^\rg} \bs\mu_{i;\dimv_i},
$$
in which the factors evidently commute.
Recall the seed
$$
\quiver_{\overline a: j \to i} = \quiver_{GL_{\dimv_i},0} \sqcup \hc{e_{\overline{a}}} \sqcup \quiver_{GL_{\dimv_j},0}.
$$
An easy computation carried out in~\cite{SS25} shows that
\begin{align*}
\bs\mu_{i;\dimv_i}(\quiver_{\overline a: j \to i}) &= -\quiver_{GL_{\dimv_i},0} \sqcup \hc{e_{\overline{a}}+\sum_{m=1}^{\dimv_i}(e_{s_{i,m}}+e_{t_{i,m}})} \sqcup \quiver_{GL_{\dimv_j},0}, \\
\bs\mu_{j;\dimv_j}(\quiver_{\overline a: j \to i}) &= \quiver_{GL_{\dimv_i},0} \sqcup \hc{e_{\overline{a}}+\sum_{m=1}^{\dimv_j}(e_{s_{j,m}}+e_{t_{j,m}})} \sqcup -\quiver_{GL_{\dimv_j},0}.
\end{align*}
Recalling~\eqref{eq:a-bar} we conclude that the Donaldson--Thomas transformation $\bs\mu_\Gamma$ for the seed $\quiver_\Gamma$ can be written as
\beq
\label{eq:DT}
\bs\mu_\Gamma = \bs\mu_{\Gamma_0} \circ \bs\mu_{\Gamma_1}.
\eeq

We finish this section with the following speculative application of the transformation~$\bs\mu_\Gamma$ to the $4d$ $\mathcal N=2$ super-symmetric gauge theories. Following~\cite{ACC+13, ACC+14}, we recall that the spectrum of the so-called BPS particles in the theory can be alalyzed with the help of the associated BPS quiver. We also recall that a chamber in the space of stability conditions which has a finite BPS spectrum yields a Donaldson--Thomas transformation of the BPS quiver, consisting of a finite chain of sign-coherent mutations. To the best of our understanding, at this point there is no known algorithm for computing a BPS quiver of a general $4d$ $\mathcal N=2$ super-symmetric quiver gauge theory. However, based on the observation that for all examples computed in~\cite{ACC+13, ACC+14} the BPS quiver coincides with the cluster quiver $\quiver_\Gamma$ and the existences of the Donaldson--Thomas transformation $\bs\mu_\Gamma$, we make the following conjecture.

\begin{conjecture}
The cluster quiver $\quiver_\Gamma$ is a BPS quiver of the $4d$ $\mathcal N=2$ super-symmetric gauge theory defined by the quiver $\Gamma$. Moreover, if $\Gamma$ does not contain 1-cycles there exists a chamber in the space of stability conditions where BPS spectrum of the corresponding theory is finite.
\end{conjecture}

\section{Cluster ensembles for unframed simply connected quivers}
\label{sec:avars}
In this section we consider the case of quiver gauge theories where the underlying graph $\Gamma$ has no framing nodes, and satisfies $H_1(\Gamma,\Z)=0$.  In this case the skew-bilinear form on $\Lambda_{|\Gamma|}$ is non-degenerate and unimodular, so we can extend the cluster seed $\quiver_\Gamma$ to a compatible pair. Since taking disjoint unions of quivers amounts to taking tensor products over $\Z[q^{\pm1}]$ of the Coulomb branch rings, it suffices to explain how to do this in the case that $\Gamma$ is a tree.

As in the proof of Lemma~\ref{lem:linalg} we choose a distinguished root vertex $i_\bullet\in \Gamma_0$, and write $\Gamma^\dag = (\Gamma,i_\bullet)$ for the corresponding rooted tree. Recall that we get in this way a basis $\{e_\ell\}$ for $\Lambda_{|\Gamma|}$ given by adjoining to $\quiver_\Gamma$ the vectors
$$
e_i = x_{i,1}, \quad i\in \Gamma_0 \qquad \text{and } \quad e_{\bullet} = p_{i_\bullet,\dimv_{i_\bullet}}.
$$
For each node $i\in\Gamma_0$, write $\Gamma^\dag_{> i}$ for the set of all descendants of $i$ in the rooted tree $\Gamma^\dag$, and set $\Gamma^\dag_{\geq i}=\Gamma^\dag_{>i}\cup\{i\}$. Given a node $i\in\Gamma_0$, set $\sgn(i)=1$ if $i$ is the root, or $i$ receives an arrow from its parent in the directed graph $\Gamma$, and set $\sgn(i)=-1$ otherwise. Similarly, if $a$ is a directed edge in $\Gamma$, define $\sgn(a)=1$ if $a$ points from parent to child, and $\sgn(a)=-1$ otherwise. Now introduce
$$
\rho(i) = \begin{cases} 
\sum_{a:i\rightarrow j}\sum_{k\in\Gamma^\dag_{\geq j}}x_{\varpi_{k,\dimv_k}}, &\sgn(i)=1, \\
-{x}_{\varpi_{i,\dimv_i}}-\sum_{a:j\rightarrow i}\sum_{k\in\Gamma^\dag_{\geq j}}x_{\varpi_{k,\dimv_k}}, &\sgn(i)=-1.
\end{cases}
$$

\begin{prop}
\label{prop:avars}
Let $\quiver_\Gamma^\vee = \hc{\xi_\bullet, \xi_a, \xi_i, \xi_{s_{i,m}}, \xi_{t_{i,m}} \,\big|\, a \in \Gamma_1, \, i \in \Gamma_0^\rg, \, 1 \le m < \dimv_i}$, where
\begin{align*}
&\xi_{t_{i,m}} = x_{\varpi_{i,m}} +\rho(i), &
&\xi_{s_{i,m}} = p_{\varpi_{i,m}} + x_{\varpi_{i,m}} - p_{\varpi_{i,\dimv_i}} +\rho(i), \\
&\xi_\bullet = \sum_{i\in\Gamma_0}x_{\varpi_{i,\dimv_i}}, &
&\xi_i = -p_{\varpi_{i,d_i}},
\end{align*}
and
$$
\xi_{a:i\to j} =
\begin{cases} 
\sum_{k\in\Gamma^\dag_{\geq j}}x_{\varpi_{k,\dimv_k}}, &\sgn(a)=1, \\
-\sum_{k\in\Gamma^\dag_{\geq i}}x_{\varpi_{k,\dimv_k}}, &\sgn(a)=-1.
\end{cases}
$$
Then $(\quiver_\Gamma, \quiver_{\Gamma}^\vee)$ is a compatible pair in the lattice $\Lambda_{|\Gamma|}$.
\end{prop}

\begin{proof}
Let us compute the elements $\xi_{t_{k,m}}$ in the basis  $\quiver_\Gamma^\vee$ dual to $\quiver_\Gamma$, with all other elements in the dual basis $\quiver_\Gamma^\vee$ being treated in a similar way. Consider the ansatz
$$
\xi_{t_{k,m}}(\beta) = x_{\varpi_{k,m}}+\sum_{r\in \Gamma_0}\beta_r x_{\varpi_{r,\dimv_r}}, \quad \beta_r \in\mathbb{Z}.
$$
Since all the $e_{s_{i,m}}$ and $e_{t_{i,m}}$ lie in the direct sum of two copies of the root lattice and thus pair to zero with all $p_{\varpi_{j,\dimv_j}}$, $x_{\varpi_{j,\dimv_j}}$, we are guaranteed that the resulting vector $\xi_{t_{k,m}}(\beta)$ has the desired pairing with all $e_{t_{i,n}}$ and is orthogonal to all vectors $e_i$ and $e_{s_{i,n}}$. We now determine the coefficients $\beta_i$ using the conditions of orthogonality with the basis vectors $e_a$ and $e_\bullet$. Taking the pairing with $e_{a:i\rightarrow j}$ gives a relation
\begin{align}
\label{eq:de}
\beta_j-\beta_i = \delta_{i,k}
\end{align}
which says that the function $\beta$ is constant along arrows unless the source of the arrow is $k$. The orthogonality with $e_\bullet$ fixes the initial condition $\beta_{i_\bullet}=0$, and since we assume $\Gamma$ connected the relation~\eqref{eq:de} then determines all other $\beta_r$. This way we arrive at the claimed formula $\sum_{r\in \Gamma_0}\beta_r x_{\varpi_{r,\dimv_r}}=\rho(i)$.
\end{proof}

\begin{cor}
\label{cor:avars}
If $\Gamma$ is a quiver without cycles or framing nodes, the quantum cluster $\Ac$-variables defined by the compatible pair $(\quiver_\Gamma,\quiver^\vee_\Gamma)$ consist of dressed minuscule monopole operators in $\Acr_\Gamma$.
\end{cor}

\begin{figure}
\begin{tikzpicture}[thick, x=.7cm, y=.9cm, every node/.style={circle, draw, fill=white, minimum size=6mm, inner sep=0pt, font=\footnotesize}]

  \node (Asi)  at (-1.8,1.2) {${s_{i,1}}$};
  \node (Ati)  at ( 1.2,1.2) {${t_{i,1}}$};
  \node (Aij)  at (-.2,0) {$a$};
  \node (Asj)  at (-1.8,-1.4) {$s_{j,1}$};
  \node (Atj)  at ( 1.2,-1.4) {$t_{j,1}$};

  \node[rectangle] (Ab)   at (-.2,2.4) {$\bullet$};
  \node[rectangle] (Ai)   at ( 3.4,1.2) {$i$};
  \node[rectangle] (Aj)   at ( 3.4,-1.4) {$j$};

  \draw[->] (Ab)  -- (Asi);
  \draw[->] (Ati) -- (Ab);
  \draw[->] (Ati) -- (Ai);
  \draw[->] (Ai)  -- (Aij);
  \draw[->] (Atj) -- (Aj);

  \draw[->] (Asi.10) to (Ati.170);
  \draw[->] (Asi.-10) to (Ati.-170);

  \draw[->] (Asj.10) to (Atj.170);
  \draw[->] (Asj.-10) to (Atj.-170);

  \draw[->] (Ati) -- (Aij);
  \draw[->] (Aij) to (Asi);

  \draw[->] (Atj) -- (Aij);
  \draw[->] (Aij) -- (Asj);
\end{tikzpicture}
\caption{Cluster quiver for the extended seed $\quiver_\Gamma$ in Example~\ref{eg:seed}.}
\label{fig:exq}
\end{figure}

\begin{example}
\label{eg:seed}
Let us spell out the $\Ac$-variables for the case $\Gamma = i\rightarrow j$, with $d_i=d_j=2$ and the vertex $i$ taken as the root. The cluster-algebra quiver representing the extended initial seed $\quiver_\Gamma$ is shown in Figure~\ref{fig:exq}, and the compatible pair reads
\begin{align*}
\Pi_\Gamma &= \hc{e_{s_{i,1}}, e_{t_{i,1}}, e_i, e_\bullet, e_{s_{j,1}}, e_{t_{j,1}}, e_j, e_a}, \\
\Pi^\vee_\Gamma &= \hc{\xi_{s_{i,1}}, \xi_{t_{i,1}}, \xi_i, \xi_\bullet, \xi_{s_{j,1}}, \xi_{t_{j,1}}, \xi_j, \xi_a},
\end{align*}
where
\begin{align*}
e_{s_{i,1}} &= -x_{\alpha_{i,1}}, & e_{t_{i,1}} &= p_{\alpha_{i,1}}+x_{\alpha_{i,1}}, & e_i &= x_{i,1}, & e_\bullet &= p_{i,2}, \\
e_{s_{j,1}} &= -x_{\alpha_{j,1}}, & e_{t_{j,1}} &= p_{\alpha_{j,1}}+x_{\alpha_{j,1}}, & e_j &= x_{j,1}, & e_a &= p_{j,2}-p_{i,1},
\end{align*}
and
\begin{align*}
\xi_{s_{i,1}} &= p_{\varpi_{i,1}^*}+x_{\varpi_{i,1}+\varpi_{j,2}}, & \xi_{t_{i,1}} &= x_{\varpi_{i,1}+\varpi_{j,2}}, & \xi_i &= -p_{\varpi_{i,2}}, & \xi_\bullet &= x_{\varpi_{i,2}+\varpi_{j,2}}, \\
\xi_{s_{j,1}} &= p_{\varpi_{j,1}^*}+x_{\varpi_{j,1}}, & \xi_{t_{j,1}} &= x_{\varpi_{j,1}}, & \xi_j &= -p_{\varpi_{j,2}}, & \xi_a &= x_{\varpi_{j,2}}.
\end{align*}
Hence we have five mutable $\Ac$-variables
\begin{align*}
Y_{\xi_{s_{i,1}}} &=  q^{-1}[\Oc_{\Rcr_{\varpi_{i,1}+\varpi_{j,2}}}\otimes \Sc^\vee_{\varpi_{i,1}}], &
Y_{\xi_{t_{i,1}}} &=  q^{-1}[\Oc_{\Rcr_{\varpi_{i,1}+\varpi_{j,2}}}], \\
Y_{\xi_{s_{j,1}}} &=  [\Oc_{\Rcr_{\varpi_{j,1}}}\otimes \Sc^\vee_{\varpi_{j,1}}], &
Y_{\xi_{t_{j,1}}} &= [\Oc_{\Rcr_{\varpi_{j,1}}}], &
Y_{\xi_a} &= q^{-1}[\Oc_{\Rcr_{\varpi_{j,2}}}].
\end{align*}
together with three frozen $\Ac$-variables
$$
Y_{\xi_i} = \hs{\mathrm{det}^\vee_i}, \qquad
Y_{\xi_j} =\hs{\mathrm{det}^\vee_j}, \qquad
Y_{\xi_\bullet} = q^{-2}[\Oc_{\Rcr_{\varpi_{i,2}+\varpi_{j,2}}}].
$$
Here $\mathrm{det}^\vee_i$ refers to the trivial line bundle on $\Rcr_{0}$ with equivariant structure given by the dual of the determinant character of $GL(V_i)$, so that $Y_{\xi_i},Y_{\xi_j}\in K_{G_\Gamma}(\pt)$.
\end{example}

\bibliographystyle{alpha}

\end{document}